\documentclass{amsart}

\usepackage{latexsym}
\usepackage{amssymb}
\usepackage{graphicx}

\newtheorem{THM}{Theorem}
\newtheorem{Lemma}[THM]{Lemma}
\newtheorem{Conj}[THM]{Conjecture}

\renewcommand{\Re}{\mathbb{R}}
\renewcommand{\P}{\mathbb{P}}
\newcommand{\be}{\begin{equation}}
\newcommand{\ee}{\end{equation}}

\renewcommand{\S}{\mathbb{S}}

\newcommand{\ol}[1]{\overline{#1}}
\newcommand{\ra}{\rightarrow}

\newcommand{\ip}[2]{\left\langle#1,#2\right\rangle}
\newcommand{\lp}{\left(}
\newcommand{\rp}{\right)}
\newcommand{\lb}{\left[}
\newcommand{\rb}{\right]}
\newcommand{\lc}{\left\{}
\newcommand{\rc}{\right\}}
\newcommand{\lab}{\left|}
\newcommand{\rab}{\right|}
\newcommand{\Lap}{\Delta}

\DeclareMathOperator{\dist}{dist}

\newcommand{\sO}{\mathcal{O}}

\newcommand{\sB}{\mathcal{B}}

\newcommand{\E}{\mathbb{E}}

\renewcommand{\d}[1]{\,d#1}
\renewcommand{\phi}{\varphi}

\newcommand{\QV}[1]{\left\langle#1\right\rangle}
\newcommand{\CV}[2]{\left\langle#1,#2\right\rangle}

\begin{document}

\title{A martingale approach to minimal surfaces}
\author{Robert W.\ Neel}
\address{Department of Mathematics, Columbia University, New York, NY}
\begin{abstract}
We provide a probabilistic approach to studying minimal surfaces in $\Re^3$.  After a discussion of the basic relationship between Brownian motion on a surface and minimality of the surface, we introduce a way of coupling Brownian motions on two minimal surfaces.  This coupling is then used to study two classes of results in minimal surface theory, maximum principle-type results, such as weak and strong halfspace theorems and the maximum principle at infinity, and Liouville theorems.
\end{abstract}
\thanks{The author gratefully acknowledges support from a Clay Liftoff Fellowship and an NSF Postdoctoral Research Fellowship.}
\email{neel@math.columbia.edu}
\date{June 18, 2008}

\maketitle

\section{Introduction}

The past several years have seen sustained interest in the theory of minimal surfaces in $\Re^3$.  We present an approach to studying minimal surfaces using Brownian motion and the methods of martingale theory.  We begin with a discussion of the basic relationship between Brownian motion on a minimal surface and the coordinate functions and the Gauss map, particularly in the cases when the surface is either stochastically complete with bounded curvature or (geodesically) complete and properly immersed.  We then introduce a way of coupling Brownian motions on a pair of minimal surfaces such that the particles are ``encouraged'' to couple in finite time.  We apply this coupling to two classes of results for minimal surfaces, maximum principle-type theorems, by which we mean weak and strong halfspace theorems and the maximum principle at infinity, and Liouville theorems.

In both cases, we discuss the relationship between the results we obtain from the coupling and the existing results and conjectures obtained by non-stochastic methods of geometric analysis.  Further results (beyond what we have been able to prove here) using the coupling appear to be possible, subject to obtaining a better understanding of the behavior of the process when the particles are close.  The reader who is primarily interested in the geometric consequences of the coupling is encouraged to proceed to Section~\ref{Section:Last}, where the topics just mentioned are discussed in detail.

In addition to the particular results we are able to prove, there are a few more general reasons for introducing this approach to minimal surfaces, aside from the subjective claim that probabilistic methods are intuitively appealing.  First, it provides a common framework for studying maximum principle-type theorems and Liouville theorems.  Second, there is reason to think that such a probabilistic framework might make it easier to extend results from minimal surfaces to minimal surfaces-with-boundary.  This is because Brownian motion can simply be stopped at the boundary, and prior to hitting the boundary it is governed by the local geometry just as in the boundary-less case.  For one example of such an extension, see Theorems~\ref{THM:WeakWeakHalfspace} and~\ref{THM:MaxAtInfinity} and their proofs.

I would like to thank Dan Stroock for numerous helpful discussions about Brownian motion and geometry and for comments on an earlier draft of this paper.  I would also like to thank Ioannis Karatzas for advice on solving martingale problems (in particular, the proof of Theorem 6.1 of~\cite{Karatzas} provided an outline for the proof of Theorem~\ref{THM:Main} below) and Michel \'Emery for comments on an earlier draft of this paper.  Finally, I am grateful to Bob Finn for introducing me to minimal surfaces several years ago.

\section{Brownian motion on a minimal surface}

\subsection{Basic results}\label{Section:Basic}

We begin by considering Brownian motion on a connected minimal surface $M$.  This is interesting in its own right, and it will also be important for understanding coupled Brownian motion later on.  There are several ways to think of Brownian motion on a manifold.  For our purposes, we will think of it as the solution to the martingale problem corresponding to one-half the Laplacian on $M$; we now explain what this means in more detail.   In general, Brownian motion may only be defined until an explosion time, which we denote by $e$.  Let $(C[0,e(\omega)), \sB)$ be the space of continuous paths on $M$, which we allow to escape to infinity at some time $e(\omega)\in(0,\infty]$ depending on the path $\omega$, with the Borel $\sigma$-algebra; we give the space the topology of uniform convergence on compacts.  Let $\{\sB_{t\wedge e}, t\geq 0\}$ be the filtration generated by these continuous paths.  Then Brownian motion started at $x\in M$ is a probability measure on $(C[0,e(\omega)), \sB)$, which we denote $P_x$, such that
\[\begin{split}
&P_{x}\lb \omega(0)=x\rb =1 \\
\text{and}\quad&\lp h(\omega(t\wedge e))-\int_0^{t\wedge e} \frac{1}{2}\Lap h(\omega(s))\, ds ; \sB_{t\wedge e} \rp\quad \text{is a $P_{x}$-martingale,}
\end{split}\]
for any smooth, compactly supported function $h$ on $M$.  Throughout this paper, all of our martingales will be continuous, and having established this convention we will refer to them simply as ``martingales.''

Minimal surfaces are characterized by the fact that the restrictions of the coordinate functions to the surface are harmonic.  From a probabilistic point of view, this means that the compositions of the coordinate functions with Brownian motion are local martingales.  More concretely, consider any Euclidean coordinate system $(x_1,x_2,x_3)$ on $\Re^3$.  The gradient of any coordinate $x_i$ (in $M$) is the projection of the gradient of $x_i$ in $\Re^3$ onto the tangent plane of $M$.  If we denote the unit normal to $M$ by $m=(m_1,m_2,m_3)$, then it is easy to check that $\ip{\nabla x_i}{\nabla x_j}=\delta_{ij}-m_i m_j$ where $\delta_{ij}$ is the Kronecker delta function.  Let $B_t$ be Brownian motion on $M$.  We will adopt the usual convention of writing the process obtained by composing a function with Brownian motion as $f\circ B_t= f_t$, and when there is no possibility of confusion we will sometimes omit the subscript.  We see that that the quadratic variations and cross variations of the coordinate processes are given by
\begin{equation}\label{Eqn:QV}
d\CV{x_{i,t}}{x_{j,t}} = \lp \delta_{ij}- m_i m_j \rp\d{t} .
\end{equation}
This determines the evolution of the coordinate processes in terms of the normal vector.  Another consequence of this equation is that the coordinate processes $x_{i,t\wedge e}$ are true (as opposed to local) martingales.

Next, we wish to understand the evolution of the normal vector under Brownian motion on $M$.  For any smooth surface in $\Re^3$, we can define the Gauss map $G:M\ra \S^2$ by identifying the normal vector at every point of $M$ with the corresponding point in $\S^2$.  It is well known that, if $M$ is minimal, the Gauss map is conformal with conformal factor $-\sqrt{|K|}$, where $K$ is the Gauss curvature.  Recall that a map is conformal with conformal factor $c$ (where $c$ can vary over the domain) if its differential preserves angles and scales lengths by a factor of $c$, where $c<0$ means that the map is orientation-reversing.  On a minimal surface $K\leq 0$, so the Gauss map is orientation-reversing and distorts area by a factor of $|K|$.  This implies that $G_t$ is time-changed Brownian motion on $\S^2$, with the changed time given by $\int_0^t |K_s| \d{s}$.

The introduction of the Gauss map raises the question of orientability of our minimal surfaces.  We do not wish to restrict ourselves to orientable surfaces.  To accomplish this, observe that the Gauss sphere process is well-defined (given a choice of normal vector at the starting point) whether or not the surface is orientable, because given any path there is a unique continuous choice of normal along it.  This allows us to state Theorem~\ref{THM:GaussProcess} below, for example, without regard to orientability.  From the geometric point of view, this corresponds to the fact that Brownian motion lifts to the orientation cover of the surface.  With this in mind, we freely ignore questions of orientability in what follows.

In order to say more, we need to place more restrictions on our minimal surface.  A natural assumption is that $M$ is (geodesically) complete.  An example of a theorem which can be proved in this generality is a result of Osserman that the Gauss map of a non-planar complete minimal surface is dense in $\S^2$ (see Theorem 8.1 of~\cite{Osserman}).  However, as mentioned above, we are more interested in global control of the immersion and on Liouville properties of the surface.  If we place no restrictions on $M$ beyond completeness, then not much can be said in this direction.  In~\cite{Nadir}, Nadirashvili showed that there exists a complete, minimal, conformal immersion of the disk into the unit ball.

Another natural assumption is that of stochastic completeness.  Probabilistically, this means that Brownian motion almost surely exists for all time, that is, it does not explode by ``going off to infinity'' in finite time.  Analytically, stochastic completeness means that the Cauchy initial value problem for one-half the Laplacian with bounded initial data has a unique bounded solution for all time.  Let $\rho=\sqrt{x_1^2+x_2^2+x_3^2}$.  Then a straightforward computation shows (see section 5.2.2 of~\cite{StroockGeo}) that
$\frac{1}{2}\Lap \rho^2 = 2$.  It follows from the defining property of the martingale problem that
\begin{equation}\label{Eqn:StoComp}
\E\lb \rho_{t\wedge e}^2\rb=\rho_0^2+2\lp t\wedge e\rp \quad\text{for all $t\in[0,\infty)$}.
\end{equation} 
If $M$ is stochastically complete (which is equivalent to $e$ being identically infinite), this equation holds with $t\wedge e$ replaced by $t$.  This obviously prevents $M$ from being contained in a ball.  Moreover, because the quadratic variation of a single coordinate is no greater than $t$, it also implies that $M$ cannot be contained in an infinite cylinder (that is, a set of the form $x_1^2+x_2^2\leq C$ for some $C>0$ and some choice of orthonormal coordinates on $\Re^3$).  In other words, if $M$ is a stochastically complete minimal surface, then there is at most one element of $\Re\P^2$ such that the projection of $M$ onto that line is not the entire line.

If we assume only that $M$ is stochastically complete (or even complete and stochastically complete), the preceding is as much as we can say in this direction.  Jorge and Xavier~\cite{JorgeXavier} show how to construct complete, minimal, conformal immersions of the disk into the ``slab'' $\{(x_1,x_2,x_3):-1\leq x_1\leq 1 \}$.  It is relatively easy to see that their construction can be performed in such a way that the resulting minimal surface is also stochastically complete.  Thus the above result that, for a stochastically complete minimal surface $M$, there is at most one line in $\Re\P^2$ such that the projection of $M$ onto that line is not the entire line is sharp.

\subsection{Bounded curvature and the weak halfspace theorem}

In~\cite{Xavier}, Xavier proved the weak halfspace theorem for complete minimal surfaces of bounded curvature; that is, any complete minimal surface of bounded curvature which is not a plane is not contained in any halfspace.  Our goal here is to prove a differential version of the weak halfspace theorem, namely, that the Gauss process accumulates infinite occupation time in every open subset of the Gauss sphere, almost surely.  Xavier's weak halfspace theorem as just mentioned is an obvious corollary, given the relationship between the normal vector and the quadratic variation of the coordinate processes in Equation~\eqref{Eqn:QV}.

It is well known that any complete manifold with bounded curvature is stochastically complete, and thus the assumptions in the theorem below are weaker than in Xavier's weak halfspace theorem.

\begin{THM}\label{THM:GaussProcess}
Let $M$ be a minimal surface, and assume either that $M$ is recurrent or that $M$ is stochastically complete and has bounded curvature.  If $M$ is not flat (that is, $K$ is not identically zero), then the corresponding Gauss sphere process almost surely accumulates infinite occupation time in each open subset of $\S^2$.
\end{THM}

\emph{Proof:} First, we observe that recurrent minimal surfaces are easy to handle.  Any recurrent manifold is necessarily stochastically complete.  Now suppose that $M$ is a recurrent, non-flat minimal surface.  Because $M$ is not flat, there are at least two points in $\S^2$ in the image of the Gauss map.  By recurrence, the normal vector visits a neighborhood of each of these points infinitely often, and thus we see that the Gauss sphere process has infinite quadratic variation.  Since it is time-changed Brownian motion, the Gauss process visits every open subset of $\S^2$ infinitely often.  For any open subset of $\S^2$, we can choose some component of its pre-image in $M$.  By recurrence, Brownian motion on $M$ accumulates infinite occupation time in that pre-image almost surely, and we conclude that the Gauss sphere process accumulates infinite occupation time in every open set of $\S^2$ almost surely.

Now assume that $M$ is transient, not flat, and stochastically complete with bounded curvature.  Thus its universal cover $\tilde{M}$ is conformally equivalent to the unit disk, by the uniformization theorem.  It is obviously enough to prove the result for $\tilde{M}$.  We will write the metric with respect to the usual Cartesian coordinates on the disk as $\lambda\delta_{ij}$, where $\delta_{ij}$ is the Kronecker delta function (so lengths are scaled by $\sqrt{\lambda}$ while area is scaled by $\lambda$).  Here, of course, $\lambda$ is smooth and positive on the open unit disk $D$; also, $\lambda$ determines the time-change taking Euclidean Brownian motion on the disk to Brownian motion on $\tilde{M}$.  In particular, the stochastic completeness of $\tilde{M}$ is equivalent to the statement that the integral of $\lambda$ along paths of the Euclidean Brownian motion (until the first hitting time of the boundary of the disk) is almost surely infinite.

Let $f$ and $g$ be the Weierstrass data for $\tilde{M}$, as described in Lemmas 8.1 and 8.2 of~\cite{Osserman}.  In particular, $g$ is the stereographic projection of the Gauss map and thus meromorphic, and $f$ is holomorphic.  Further, we have that (see Equation (8.7) and Lemma 9.1 of~\cite{Osserman}, and note that our definition of $\lambda$ differs from Osserman's by a power of two)
\[
\lambda = \lb\frac{|f|\lp 1+|g|^2\rp}{2}\rb^2 \quad\text{and}\quad
-K = \lb\frac{4|g'|}{|f|\lp 1+|g|^2\rp^2}\rb^2 .
\]
We wish to show that the integral of $-K$ along Brownian paths of $\tilde{M}$ is almost surely infinite.  Note that $g_t$ (which, we recall, is $g\circ B_t$ where $B_t$ is Brownian motion on $\tilde{M}$) is a complex martingale with quadratic variation given by the integral of $-K$ along $B_t$.  Thus, showing that $-K$ has infinite integral along $B_t$ is equivalent to showing that $g_t$ almost surely does not converge, which in turn is the same as showing that $g$ does not have non-tangential limits at the boundary of the disk, except possibly on a set of Lebesgue measure zero.  This last equivalence is a consequence of Doob's form of the Fatou theorem.

We proceed by contradiction, assuming that with probability $2p>0$, $g_t$ does converge.  We may assume that our Brownian motion begins at the center of the disk, and thus that the hitting measure of the boundary is proportional to its Lebesgue measure.  Possibly after a rotation (of our coordinates on $\Re^3$), we can assume that $g_t$ has limits with absolute value less than $C$ with probability $p$, for some positive constant $C$.  Because of the above formula for $-K$, bounded curvature implies that there is a constant $C'$ such that $|g'/f|< C'$ whenever $|g|< 2C$.  Because a meromorphic function on the disk has a non-tangential limit at a boundary point if and only if it is non-tangentially bounded at that point (again by Doob's version of the Fatou theorem), it follows that $g'/f$ has a non-tangential limit with absolute value less than $C'$ on a set of boundary points of Lebesgue measure $2\pi p$.  The set of such points where the limit is zero must have Lebesgue measure zero since otherwise $g'/f$ would be identically zero, which would mean that $g'$ was identically zero, contradicting the assumption that $\tilde{M}$ is not flat.  So there is a set $\Phi$ of boundary points of Lebesgue measure $2\pi p$ where $g$ has non-tangential limits of absolute value less than $C$ and $g'/f$ has non-tangential limits with absolute value in $(0,C')$.

For any Euclidean Brownian path hitting the boundary at $\Phi$, we see that the integral of $\lambda$ along the path being infinite almost surely implies that the integral of $|f|^2$ along the path is infinite, using the above formula for $\lambda$ in terms of the Weierstrass data.  This in turn almost surely implies that the integral of $|g'|^2$ along the path is infinite.  However, the integral of $|g'|^2$ is the quadratic variation of the $g_t$ (which, we recall, is a complex martingale).  We conclude that $g_t$ almost surely does not converge along Brownian paths which hit the boundary at $\Phi$.  This is a contradiction, and thus we have shown that the integral of $-K$ along Brownian paths of $\tilde{M}$ is almost surely infinite.

Because the Gauss sphere process is time-changed Brownian motion on the sphere with the time-change given by the integral of $-K$, this shows that the Gauss sphere process hits every open subset of $\S^2$ infinitely often.  Now consider $B_{\epsilon}$ an open ball in $\S^2$ of radius $\epsilon$, and let $B_{\epsilon/2}$ be an open ball with the same center and half the radius.  Because $-K$ is bounded from above, there exists some $\delta>0$ such that every time the Gauss sphere process hits $B_{\epsilon/2}$, it spends time $\delta$ in $B_{\epsilon}$ with probability at least $1/2$.  An easy application of the Borel-Cantelli lemma then shows that the process spends an infinite amount of time in $B_{\epsilon}$, almost surely.  As this argument applies to any open ball in $\S^2$, the theorem is proved. $\Box$

We note that this theorem does not require that $M$ is complete.  This also explains why the hypothesis of the theorem is that $M$ is not flat, rather than that $M$ is not a plane.  For example, let $M$ be the universal cover of a plane minus two points, with the obvious immersion.  Then $M$ is a flat, stochastically complete minimal surface, but $M$ is transient and, in particular, not a plane.

\subsection{Proper immersion and the weak halfspace theorem}

A weak halfspace theorem can also be proved for properly immersed minimal surfaces; namely, if $M$ is a complete, properly immersed minimal surface, and $M$ is not a plane, then $M$ is not contained in any halfspace.  This was first done by Hoffman and Meeks~\cite{HoffmanMeeks}, using a geometric construction comparing $M$ to the lower half of a catenoid.  As noted in~\cite{MeeksSurvey} (see the discussion surrounding Theorems 1.3 and 1.4 there), it is also a simple consequence of Theorem 3.1 of~\cite{CKMR}, which is proved using elementary harmonic function methods.

Before continuing, we make an observation about terminology.  Any properly immersed minimal surface is necessarily complete, and thus the phrase ``complete, properly immersed minimal surface'' is somewhat redundant.  Nonetheless, we will employ this phrase in order to highlight the parallel between them and stochastically complete minimal surfaces of bounded curvature, since we are viewing both conditions as strengthenings of the corresponding assumptions of (geodesic or stochastic) completeness.

We claim that Equation~\eqref{Eqn:StoComp} implies that any complete, properly immersed minimal surface is stochastically complete.  Brownian motion on a complete manifold explodes in finite time if and only if it exits every compact set in finite time.  Properness of the immersion means that exiting every compact set of $M$ is the same as exiting every compact set of $\Re^3$, which is the same as $\rho^2_t$ blowing up in finite time.  However, Equation~\eqref{Eqn:StoComp} implies that the expectation of $\rho^2_t$ remains finite at all times.  This justifies our claim.  Thus the assumption that $M$ is properly immersed (and, necessarily, complete) is a strengthening of the assumption that $M$ is stochastically complete.

Our probabilistic proof is similar, both in spirit and technique, to the harmonic function approach mentioned above.  Also, we note that we are unable to give a stronger, ``differential,'' version as we did for the bounded curvature case (more on this below).

\begin{THM}
Let $M$ be a complete, properly immersed minimal surface.  If $M$ is not flat, then $M$ is not contained in any halfspace.
\end{THM}

\emph{Proof:}  Assume that $M$ is not flat.  The case when $M$ is recurrent is already covered by Theorem~\ref{THM:GaussProcess}, so we assume that $M$ is transient.

We begin by showing that the integral of $1-m_3^2$ along Brownian motion on $M$ blows up almost surely.  We will proceed by contradiction; assume that the integral of $1-m_3^2$ is bounded with probability $p>0$.  Let $r=\sqrt{x_1^2+x_2^2}$.  Then the semi-martingale decomposition of the $r_t$ process is
\[
dr_t = \sqrt{1-\lp \partial_r \cdot m\rp^2} \d{W_t} +\frac{1}{2r_t} \lp m_3^2+\lp \partial_r \cdot m\rp^2 \rp \d{t} +\d{L_t}
\]
where $W_t$ is some Brownian motion, and $L_t$ is an increasing process which increases only when $r_t=0$.  Also, $\partial_r$ is the $\Re^3$-gradient of $r$ (since we can, at least locally, identify $M$ with its image in $\Re^3$, there is no problem in thinking of Brownian motion on $M$ as a martingale in $\Re^3$ and writing the coefficients of the semi-martingale decomposition in terms of the geometry of $\Re^3$).  This is the geometer's convention of identifying vector fields with first-order differential operators, and we will adopt this convention freely in what follows.  We will assume, for simplicity, that we start our Brownian motion at a point with $r_0>0$.

We wish to compare the $r_t$ process with the 2-dimensional Bessel process generated by the same Brownian motion, that is, the strong solution of
\[
\d{\tilde{r}_t} = \d{W_t} + \frac{1}{2\tilde{r}_t}\d{t} \quad\text{with $\tilde{r}_0=r_0> 0$.}
\]
Then the process $(r-\tilde{r})_t$ has the semi-martingale decomposition
\[
\d{\lp r-\tilde{r} \rp_t}= \lb  \sqrt{1-\lp \partial_r \cdot m\rp^2} -1\rb \d{B_t} 
+\lb \frac{1}{2r_t} \lp m_3^2+\lp \partial_r \cdot m\rp^2 \rp -\frac{1}{2\tilde{r}_t}\rb  \d{t} +\d{L_t}.
\]
We will also write this in integrated form as $(r-\tilde{r})_t=M_t+A_t+L_t$.
We wish to control the size of excursions of $r_t$ above $\tilde{r}_t$.  Introduce the stopping time
\[
\sigma(C)=\inf\left\{t\geq0 : \int_0^t\lp 1-m_3^2\rp \geq C\right\} .
\]
We will assume that $C$ is chosen large enough so that $\sigma(C)=\infty$ with probability at least $3p/4$.

The quadratic variation of $M$ satisfies
\[\begin{split}
\QV{M}_{t\wedge\sigma(C)} &= \int_0^{t\wedge\sigma(C)} \lb 2\lp 1-\sqrt{1-\lp \partial_r \cdot m\rp^2}\rp -\lp \partial_r \cdot m\rp^2\rb \d{\tau} \\
&\leq 2 \int_0^{t\wedge\sigma(C)} \lp 1-\sqrt{1-\lp \partial_r \cdot m\rp^2}\rp \d{\tau} .
\end{split}\]
Further, because $\lp \partial_r \cdot m\rp^2\leq 1-m_3^2$, we see that
\[
 1-\sqrt{1-\lp \partial_r \cdot m\rp^2} \leq 1-|m_3| \leq 1-m_3^2 .
\]
It follows that $\QV{M}_{t\wedge\sigma(C)}\leq C $.  Thus whatever else it is, $M_{t\wedge\sigma(C)}$ is a continuous process with both $\sup_{t\geq0}\{ M_{t\wedge\sigma(C)}\}$ and $\inf_{t\geq0}\{ M_{t\wedge\sigma(C)}\}$
almost surely finite.

Now introduce $\eta(t)$ as the random (but not stopping) time defined by
\[
\eta(t)=\sup\{\tau\leq t : r_{\tau}\leq \tilde{r}_{\tau}+1\} ;
\]
that is, $\eta(t)$ is the last time before $t$ when $r$ has been at or below $\tilde{r}+1$.
Since integrals of the non-martingale terms make sense with random times that aren't stopping times, we can write
\[\begin{split}
\lp r-\tilde{r} \rp_{t\wedge\sigma(C)}&\leq 1+ M_{ t\wedge\sigma(C)}-M_{\eta(t\wedge\sigma(C))} \\
&+\int_{\eta(t\wedge\sigma(C))}^{ t\wedge\sigma(C)} \lb \frac{1}{2r_{\tau}} \lp m_3^2+\lp \partial_r \cdot m\rp^2 \rp -\frac{1}{2\tilde{r}_{\tau}}\rb  \d{\tau} + \int_{\eta(t\wedge\sigma(C))}^{ t\wedge\sigma(C)}\d{L_{\tau}}.
\end{split}\]
Note that $m_3^2+\lp \partial_r \cdot m\rp^2 \leq 1$, and thus that the first integrand is always non-positive on the interval of integration.  Next, because $\tilde{r}_t$ is a 2-dimensional Bessel process, it is almost surely positive for all time.  Thus the $L_t$ term never increases when $r_t\geq \tilde{r}_t$, and thus never increases on the interval of integration.  It follows that
\[
\lp r-\tilde{r} \rp_{t\wedge\sigma(C)}\leq 1+ M_{ t\wedge\sigma(C)}-M_{\eta(t\wedge\sigma(C))}
\]
Then we conclude that
\[\begin{split}
\sup_{t\geq0}\{r_{t\wedge\sigma(C)} - \tilde{r}_{t\wedge\sigma(C)}\} &<1+ \sup_{t\geq0}\{ M_{t\wedge\sigma(C)}\} - \inf_{t\geq0}\{ M_{t\wedge\sigma(C)}\} \\
&< \infty
\end{split}\]
almost surely.

To continue, note that because 
\[
\QV{x_3}_{t\wedge\sigma(C)}=\int_0^{t\wedge\sigma(C)}(1-m_3^2)\d{\tau}<C ,
\]
we see that $\sup_{t\geq0}\{|x_3|_{t\wedge\sigma(C)}\}<\infty$ almost surely.
Recall that with probability $3p/4>0$, we have $\sigma(C)=\infty$. Thus with probability $3p/4$, we have that $\sup_{t\geq0}\{r_{t} - \tilde{r}_{t}\} <\infty$ and $\sup_{t\geq0}\{|x_3|_{t}\}<\infty$.  Now $\tilde{r}_t$ is recurrent for the set $\tilde{r}<1$.  It follows that, on the set of paths with $\sigma(C)=\infty$, we know that $r_t$ returns infinitely often to some bounded interval (depending on the particular path).  Since a subset of $\Re^3$ with $r$ and $x_3$ bounded is contained in a compact cylinder, it follows that the set of paths with $\sigma(C)=\infty$ returns infinitely often to a compact cylinder (the exact cylinder depending on the path).  Because the immersion is proper, the intersection of any compact subset of $\Re^3$ with $M$ is also compact.  Thus, almost every path with $\sigma(C)=\infty$ returns infinitely often to some compact subset of $M$.  This implies that, with probability $3p/4$, Brownian motion on $M$ is recurrent.  This contradicts our assumption that $M$ is transient, and we conclude that the integral of $1-m_3^2$ is almost surely unbounded.

This means that $x_{3,t}$ is  a martingale with almost surely unbounded quadratic variation.  Thus $M$ is not contained in the halfspace $x_3>0$.  Since our choice of Euclidean coordinates on $\Re^3$ was arbitrary, it follows that $M$ cannot be contained in any halfspace.  $\Box$

Note that the computation showing that Brownian motion on $M$ is recurrent on some compact subset of $\Re^3$ doesn't use the fact that the process is restricted to a surface.  In particular, consider any $\Re^3$-valued (continuous) martingale $Y$ such that the diffusion matrix of $Y$ at each instant is given by the diffusion matrix for Brownian motion on some plane (if $Y$ is Brownian motion on $M$, then this plane is the tangent plane to $M$ at $Y$).  Further, suppose that one of the coordinates of $Y$ has finite quadratic variation with probability $p>0$.  Then the proof shows that, up to a set of probability zero, those paths for which a coordinate of $Y$ is bounded are recurrent on some compact subset of $\Re^3$ (with the particular subset depending on the path).  That the corresponding plane-field is integrable (in the sense of the Frobenius theorem) and the resulting surface complete and properly immersed in the case of Brownian motion on a minimal surface is only used to show that recurrence on $\Re^3$ implies that Brownian motion on $M$ is recurrent.

One consequence of this is that the theorem also applies to complete, properly immersed branched minimal surfaces.  In particular, there are at most countably many branch points, so Brownian motion almost surely avoids them all.  We then see that the above proof applies in this case as well.

As mentioned, we do not have a ``differential'' version of this result, as we do in the bounded curvature case.  In particular, nothing in this theorem rules out the Gauss process spending only finite time in any closed set in $\Re\P^2$ not containing the North/South pole.  As we will see later, this difference will mean that we can prove more in the bounded curvature case than in the properly immersed case.

\section{Coupled Brownian motion}

To address more sophisticated issues, we will need to do more than study a single Brownian motion on a minimal surface.  In particular, for any two stochastically complete minimal surfaces $M$ and $N$ (where we allow the possibility that they are the same surface), we wish to couple Brownian motions on the surfaces such that the $\Re^3$-distance between the two particles goes to zero as efficiently as possible.  We will do so by constructing a diffusion on the product space that, pointwise, gives a favorable evolution for the distance between the particles, in a sense to be made precise below.

\subsection{Pointwise specification of an optimal coupling}

We begin by determining what the coupling should be ``instantaneously.''

Let $x=x_t$ be the position of the Brownian motion on $M$ and $y=y_t$ be the position of the Brownian motion on $N$.  Let $r=r(x,y)$ be the $\Re^3$-distance between $x$ and $y$. We assume that $r(x,y)>0$, since we stop the process if the particles successfully couple.  The instantaneous evolution of $r$ will depend only on $r$ and the positions of the tangent planes $T_xM$ and $T_yN$ (and the coupling).  We wish to understand this dependence in detail and use it to see how to choose the coupling.  As such, we will consider $x$ and $y$ to be given, fixed points and construct our analysis around them.

We can choose orthonormal coordinates $(z_1, z_2, z_3)$ for $\Re^3$ which are well-suited to the current configuration (we use $z_i$ for our coordinates here instead of $x_i$ in order to avoid overburdening the notation, considering that we use $x$ to denote a point in $M$).  Let $\partial_{z_3}$ lie in the direction of $x-y$, where we view $y$ and $x$ as elements of $\Re^3$ under the corresponding immersions of $N$ and $M$.  Let $\partial_{z_2}$ lie in $T_xM$.  (As long as $T_xM$ is not perpendicular to $x-y$, this will determine $\partial_{z_2}$ up to sign.)  Finally, $\partial_{z_1}$ is chosen to complete the orthonormal frame.  Such a choice of orthonormal frame determines orthonormal coordinates on $\Re^3$ up to translation.  Since we are only interested in the relative positions of $x$ and $y$, any choice of coordinates $\{z_i\}$ corresponding to this choice of frame $\{\partial_{z_i}\}$ will work.

Given such coordinates, $T_xM$ is determined by the angle it makes with the $z_1z_2$-plane; call this angle $\theta$.  We choose an orthonormal basis $\{\partial_{\alpha}, \partial_{\beta}\}$ for $T_xM$ such that, at $x$,
\[
\nabla_M z_1 = \cos\theta \partial_{\alpha}  , \quad
\nabla_M z_2 = \partial_{\beta}  , 
\quad\text{and}\quad\nabla_M z_3 = \sin\theta \partial_{\alpha} , 
\]
where $\nabla_M z_i$  is the gradient of $z_i$ restricted to $M$.  In other words, $\partial_{\beta}$ is in the $\partial_{x_2}$ direction, while $\partial_{\alpha}$ lies in the intersection of $T_xM$ and the $z_1z_3$-plane.

In order to specify the direction of $T_yN$, we will need two angles.  Let $\phi$ be the angle $T_yN$ makes with the $z_1z_2$-plane.  Let $\psi$ be the angle between the intersection of $T_yN$ with the $z_1z_2$-plane and $\partial_{z_2}$.  Then we can choose an orthonormal basis $\{\partial_a, \partial_b\}$ for $T_yN$ such that, at $y$,
\[\begin{split}
\nabla_N z_1 = \cos\phi\cos\psi \partial_a &+\sin\psi \partial_b ,\quad
\nabla_N z_2 = -\cos\phi\sin\psi \partial_a + \cos\psi\partial_b  , \\
&\quad \text{and}\quad\nabla_N z_3 = \sin\phi \partial_a , 
\end{split}\]
This means that $\partial_b$ lies in the intersection of $T_yN$ and the $x_1x_2$-plane, while $\partial_a$ lies along the projection of $\partial x_3$ onto $T_yN$.

Note that, after possibly reflecting some of the $z_i$ and exchanging the roles of $M$ and $N$, we can, and will, assume that
\[
\theta\in[0,\pi/2] , \quad\phi\in[0,\pi/2] , \quad\text{and } \psi\in[0,\pi] .
\]
Further, whenever $\theta$ or $\phi$ is zero, we can, and will, take $\psi$ to be zero.  With these conventions, there is a unique choice of $(\theta,\phi, \psi)$ for every $(x,y)\in (M\times N)\setminus\{r=0\}$.  Since the coupling at each point will depend only on these three angles, we will refer to them as the configuration of the system at this instant.  It is clear that the map into the configuration space is smooth near any point where all three angles are in the interior of their ranges, but in general will be discontinuous if any of them is at the boundary of its range.  For the moment, we are concerned only with determining the coupling pointwise.  Later, when we consider the induced martingale problem on $M\times N\setminus\{r=0\}$, we will have to account for the behavior of the configuration as $(x,y)=(x_t,y_t)$ varies.

With the product metric on $M\times N$, again using the fact that $x$ and $y$ can be viewed as elements of $\Re^3$, and viewing the coordinate $z_i$ as a function on $\Re^3$ that gives the $i$th coordinate, we have
\begin{equation}\label{Eqn:SystemEvo}\begin{split}
\nabla_{(M\times N)} z_1(x-y) &=\cos\theta\partial_{\alpha} -\cos\phi\cos\psi \partial_a -\sin\psi \partial_b , \\
\nabla_{(M\times N)} z_2(x-y) &= \partial_{\beta} +\cos\phi\sin\psi \partial_a - \cos\psi\partial_b  , \\
\text{and}\quad \nabla_{(M\times N)} z_3(x-y) &= \sin\theta \partial_{\alpha} -\sin\phi \partial_a . \\
\end{split}\end{equation}

Now we wish to see what the above computations mean for the evolution of $r$ under some (to be determined) coupling of Brownian motions on $M$ and $N$.  Our orthonormal basis $\{\partial_{\alpha}, \partial_{\beta}\}$ for $T_xM$ uniquely determines normal coordinates $(\alpha, \beta)$ in a neighborhood of $x$ in $M$, and we similarly have normal coordinates $(a,b)$ in a neighborhood of $y$ in $N$.  Further, $(\alpha,\beta,a,b)$ gives product normal coordinates in a neighborhood of $(x,y)$ in $M\times N$.  Since the $z_{i,t}=z_i(x_t-y_t)$ are martingales, their instantaneous evolution is determined by the rate of change of their quadratic variations.  The above expressions for the gradients mean that, at this instant with $(x_t,y_t)=(x,y)$, we have
\[\begin{split}
d\QV{z_{1,t}} &=\lp\cos^2\theta +\cos^2\phi\cos^2\psi+\sin^2\psi\rp dt -2\cos\theta
\lb \cos\phi\cos\psi\, d\CV{\alpha_t}{a_t} \right. \\
  &\quad \left. +\sin\psi\, d\CV{\alpha_t}{b_t}\rb , \\
d\QV{z_{2,t}} &= \lp1+\cos^2\phi\sin^2\psi+\cos^2\psi\rp dt -2
\lb \cos\psi\, d\CV{\beta_t}{b_t} -\cos\phi\sin\psi\, d\CV{\beta_t}{a_t}\rb , \\
\text{and}\quad &d\QV{z_{3,t}} = \lp\sin^2\theta+\sin^2\phi\rp dt -2\sin\theta\sin\phi\, d\CV{\alpha_t}{a_t} .
\end{split}\]
Here we have used that the marginals are Brownian motions, so that the only unknown quantities are the four cross-variations above, which relate the Brownian motions on $M$ and $N$.  It is these four cross-variations which we think of as determining the coupling at $(x,y)$.

The $z_i$ were chosen so that $z_1(x-y)=z_2(x-y)=0$ and $z_3(x-y)=r$.  Thus, the semi-martingale decomposition of $r=\sqrt{z_1^2+z_2^2+z_3^2}$, at this instant with $(x_t,y_t)=(x,y)$, is easy to compute using the above equations for the quadratic variations of the $z_i$.  In particular, It\^o's formula implies that, at this instant with $(x_t,y_t)=(x,y)$, the bounded variation part is growing at rate $\lp d\QV{z_{1,t}}+d\QV{z_{2,t}}\rp/2r$ while the quadratic variation of the martingale part is growing at rate $d\QV{z_{3,t}}$.  We introduce one more bit of notation.  Let $f$ and $g$ be such that the martingale part has quadratic variation given by the integral of $f$ along paths and the bounded variation part is given by the integral of $g/2r_t$ along paths (plus $r_0$).  Then $f$ and $g$ are non-negative functions of the configuration $(\theta,\phi,\psi)$ and the coupling, and they determine the evolution of $r_t$.

We return to discussing how the Brownian motions should be coupled at $(x,y)$.  We wish to consider couplings which make the cross-variations as large (in absolute value) as possible, since it is intuitively clear that such couplings give the most control over $f$ and $g$.  We are working within the framework of the martingale problem, but the type of couplings we are considering can perhaps be explained best from the point of view of stochastic differential equations.  From this perspective, we are saying that the Brownian motion on $N$ should be driven by the Brownian motion on $M$ (at least at the instant under consideration), or vice versa, since the roles of the two surfaces are symmetric.  At any rate, such couplings are, at a point, parametrized by isometries between $T_xM$ and $T_yN$.  In terms of our orthonormal bases, all such maps can be written in the form
\[\begin{split} 
\partial_a&=A\cos\sigma\partial_{\alpha} + A\sin\sigma\partial_{\beta} , \\
\partial_b&=-\sin\sigma\partial_{\alpha} + \cos\sigma\partial_{\beta} ,
\end{split}\]
where $\sigma\in[0,2\pi)$ and $A$ is $\pm 1$.  Note that we're simply relating the tangent planes by an element of $O(2)$.  Our choice of orthonormal bases for the tangent planes determines coordinates for $O(2)$ as above, where $A$ determines whether we're in the orientation-preserving or orientation-reversing component and $\sigma$ then parametrizes the relevant component, which is diffeomorphic to $\S^1$ (it is a slight extension of the term ``coordinates'' to call $\sigma$ and $A$ coordinates, but the meaning is clear and no harm is done).

In terms of $\sigma$ and $A$, the four cross-variation terms relating the Brownian motions on $M$ and $N$, at the instant when $(x_t,y_t)=(x,y)$, are given by
\begin{equation}\label{Eqn:Coupling}\begin{split}
d\CV{\alpha_t}{a_t} &= A\cos\sigma\, dt , \quad d\CV{\beta_t}{a_t}= A\sin\sigma\, dt, \\
d\CV{\alpha_t}{b_t} &= -\sin\sigma\, dt , \quad\text{and}\quad d\CV{\beta_t}{b_t}= \cos\sigma\, dt .
\end{split}\end{equation}
Then, at the instant when $(x_t,y_t)=(x,y)$, the instantaneous evolution of $r_t$ is determined by the equations
\[\begin{split}
f&= \sin^2\theta+\sin^2\phi -2A\sin\theta\sin\phi\cos\sigma \quad\text{and}\\
g&= 2+\cos^2\theta+\cos^2\phi-2\cos\sigma\lb A\cos\theta\cos\phi\cos\psi +\cos\psi\rb \\
&+2\sin\sigma\lb A\cos\phi\sin\psi+\cos\theta\sin\psi \rb .
\end{split}\]

We are now in a position to indicate what the criterion for our coupling to be ``optimal.''
We wish to consider the $r_t$ process on new a time-scale, namely the time-scale in which its quadratic variation grows at rate 1.  In this time-scale, the martingale part is just a Brownian motion, and the drift coefficient is given (at the instant for which $(x_t,y_t)=(x,y)$) by $\frac{g/f}{2r}$.  Thus, the time-changed process will be dominated by a two-dimensional Bessel process if $g\leq f$, and the domination will be strict whenever this inequality is strict.  Note that, unlike $g/f$, this inequality makes sense even if $f=0$.  That the two-dimensional Bessel process is the critical case will be seen in the following analysis (we mention it here for future reference).  That the two-dimensional Bessel process is critical is obviously significant in our effort to make $r$ hit zero as efficiently as possible, in light of the well-known relationship between the dimension of a Bessel process and its long-time behavior.

Motivated by this, we will choose our coupling so as to maximize
\begin{equation}\label{Eqn:FMinusG}\begin{split}
f-g=& -2\cos^2\theta-2\cos^2\phi +2\cos\sigma\lb A\cos\theta\cos\phi\cos\psi \right. \\
& \quad\left. +\cos\psi-A\sin\theta\sin\phi\rb-2\sin\sigma
\lb A\cos\phi\sin\psi+\cos\theta\sin\psi \rb  .
\end{split}\end{equation}
If we assume for a moment that $A$ is fixed, then it is clear how to choose $\sigma$ in order to maximize this expression.  In particular, we can view the two bracketed expressions as being the two components of a vector in the plane, in which case the two terms depending on $\sigma$ become the dot product of this vector with a unit vector making angle $\sigma$ with the positive horizontal axis.  Thus we should choose $\sigma$ so that these vectors are parallel.  Let $\sigma_{+}$ be this optimal choice when $A=1$ and $\sigma_{-}$ be this optimal choice when $A=-1$.

Writing explicit formulas for $\sigma_{+}$ and $\sigma_{-}$ isn't important (although it could easily be done in terms of inverse tangents); instead we can write directly that the maximum of $f-g$ arising from the optimal choice of $\sigma_{+}$ or $\sigma_{-}$ is
\[\begin{split}
&-2\cos^2\theta -2\cos^2\phi+2\lc \cos^2\psi \lp1+\cos^2\theta\cos^2\phi\rp +\sin^2\psi \lp\cos^2\theta+\cos^2\phi\rp \right. + \\
& \left. \sin^2\theta\sin^2\phi-2\cos\theta\cos\phi\cos\psi\sin\theta\sin\phi 
+2A\lp \cos\theta\cos\phi-\cos\psi\sin\theta\sin\phi\rp\rc^{1/2} ,
\end{split}\]
where we've used the fact that $A^2=1$.  Further, it's now clear how to choose $A$; we want $A$ to be equal to the sign of the expression it multiplies.  Doing so shows that the maximum of $f-g$, that is, the value realized by the optimal coupling, is
\[\begin{split}
&-2\cos^2\theta -2\cos^2\phi+2\lc \cos^2\psi \lp1+\cos^2\theta\cos^2\phi\rp +\sin^2\psi \lp\cos^2\theta+\cos^2\phi\rp \right. + \\
& \left. \sin^2\theta\sin^2\phi-2\cos\theta\cos\phi\cos\psi\sin\theta\sin\phi
 +2\lab \cos\theta\cos\phi-\cos\psi\sin\theta\sin\phi\rab\rc^{1/2} .
\end{split}\]
Let $\Sigma_{+}$ be the region (in $\lp M\times N\rp\setminus\{r=0\}$) where $\cos\theta\cos\phi-\cos\psi\sin\theta\sin\phi$ is non-negative, and let $\Sigma_{-}$ be the region where this expression is non-positive.  Then $\Sigma_{+}$ is the set where the maximum of $f-g$ is realized by a coupling coming from the orientation-preserving component of $O(2)$, and $\Sigma_{-}$ is the analogous set for the orientation-reversing component of $O(2)$.  Here the notions of orientation-preserving and orientation-reversing are defined relative to the standard choice of orthonormal bases for $T_xM$ and $T_yN$ given above (which we recall are not continuous when the configuration is at the boundary of its range).  Further, let $\Sigma_0$ be the set where $\cos\theta\cos\phi-\cos\psi\sin\theta\sin\phi$ is zero; this is the set on which there are two possible choices for an optimal coupling, one from each component of $O(2)$.  We will return to questions of how the choice of optimal coupling varies over $\lp M\times N\rp\setminus\{r=0\}$ below.  In particular, we will see in Lemma~\ref{Lemma:SmoothL} that $\Sigma_0$ is defined invariantly, that is, it does not depend on the particular identification of maps from $T_xM$ to $T_yN$ with $O(2)$.

Next, we wish to show that the optimal coupling is good enough for our purposes.  In particular, we wish to show that the above expression for the maximum of $f-g$ is non-negative for any values of the parameters, which will show that the time-changed $r_t$ process is dominated by a two-dimensional Bessel process (assuming that a coupling satisfying this pointwise specification exists, a question which we continue to postpone).  We also wish to determine the values of the parameters for which the expression is zero, since those will be the configurations where the process looks instantaneously like a two-dimensional Bessel process, rather than being strictly dominated by one.  Observe that the non-negativity of $f-g$ is equivalent to the inequality
\be\begin{split}\label{Eq:CouplingInequality}
\cos^2\psi & \lp1+\cos^2\theta\cos^2\phi\rp +\sin^2\psi \lp\cos^2\theta+\cos^2\phi\rp +\sin^2\theta\sin^2\phi  \\
& -2\cos\theta\cos\phi\cos\psi\sin\theta\sin\phi +2\lab \cos\theta\cos\phi-\cos\psi\sin\theta\sin\phi\rab \\
 & \quad\quad  \geq (\cos^2\theta+\cos^2\phi)^2 ,
\end{split}\ee
and positivity is equivalent to strict inequality.

We begin by considering the case when $\psi=0$ (this occurs when $x-y$ and the normal vectors to $T_xM$ and $T_yN$ are coplanar, as vectors in $\Re^3$).  Then the above simplifies to 
\[\begin{split}
& 1+\cos^2\theta\cos^2\phi +\sin^2\theta\sin^2\phi  -2\cos\theta\cos\phi\sin\theta\sin\phi +2\lab\cos(\theta+\phi)\rab \\
 & \quad\quad  \geq (\cos^2\theta+\cos^2\phi)^2 .
\end{split}\]
After simplifying the left-hand side and taking the square root of both sides, we see that this is equivalent to
\[
\lab\cos(\theta+\phi)\rab+1\geq \cos^2\theta+\cos^2\phi .
\]
One then can show that this inequality holds, and that one has equality exactly when $\theta=\phi\leq\pi/4$ or $\theta+\phi=\pi/2$.

Recall that whenever either $\theta$ or $\phi$ equals zero, we can assume that $\psi=0$.  Thus we now consider the case when all three angles are positive.  First we assume that $\theta+\phi\geq\pi/2$.  We can rewrite inequality~\eqref{Eq:CouplingInequality} as
\[\begin{split}
\cos^2\theta &+\cos^2\phi - \lp \cos^2\theta+\cos^2\phi \rp^2 
+\lp\cos^2\psi+1 \rp\sin^2\theta\sin^2\phi \\ & - 2\cos\psi\cos\theta\cos\phi\sin\theta\sin\phi 
 +2\lab \cos\theta\cos\phi - \cos\psi\sin\theta\sin\phi \rab \geq 0 .
\end{split}\]
Now $\theta+\phi\geq\pi/2$ implies that
\[
\sin\theta\sin\phi\geq\cos\theta\cos\phi\quad\text{and}\quad \cos^2\theta +\cos^2\phi \leq 1 .
\]
It follows that the sum of the first three terms is non-negative.  Along with the fact that 
\[
\cos^2\psi+1> 2\cos\psi \quad\text{for $\psi\in(0,\pi]$,}
\]
it follows that the sum of the fourth and fifth terms is strictly positive.  Since the last term is non-negative (it is an absolute value), we conclude that the inequality holds strictly.

We continue to assume that all three angles are positive, and now we consider the case $\theta+\phi<\pi/2$.  We can rewrite inequality~\eqref{Eq:CouplingInequality} as
\[\begin{split}
& 2\cos\theta\cos\phi\lp 1-\lp \cos\theta\cos\phi+\sin\theta\sin\phi\cos\psi\rp\rp 
+\cos^2\theta\sin^2\theta \\
&\quad+\cos^2\phi\sin^2\phi +\lp \cos^2\psi+1\rp\sin^2\theta\sin^2\phi
-2\sin\theta\sin\phi\cos\psi \geq 0 .
\end{split}\]
We know that this holds for $\psi=0$.  We compute the derivative of the left-hand side with respect to $\psi$ as
\[
\frac{\partial}{\partial\psi}\text{l.h.s.} = 2\sin\psi\sin\theta\sin\phi\lp 1+\cos\theta\cos\phi - \cos\psi\sin\theta\sin\phi \rp .
\]
This is positive as soon as $\psi$ is positive, and we conclude that the desired inequality holds strictly when all three angles are positive.

This shows that, under the optimal choice of $A$ and $\sigma$, $f-g$ is always non-negative.  Further, $f-g$ is positive except in the following configurations.  (That these configurations are the only ones with $f-g=0$ is something we computed above.  The additional details we are about to provide in these cases follow from elementary trigonometry and the above equations.)  When $\psi=0$ and $\theta=\phi<\pi/4$, we have $f=g=0$.  When $\psi=0$ and $\theta+\phi=\pi/2$, there are two possible couplings giving the maximum value of $f-g$, one orientation-preserving and one orientation-reversing.  For both choices of optimal coupling, $f$ and $g$ are equal and positive (except when $\theta=\phi=\pi/4$), although their shared value is larger under the orientation-reversing coupling.  Finally, the case $\psi=0$ and $\theta=\phi=\pi/4$ lies at the border between the above configurations.  Here, there are again two possibilities for the optimal coupling.  The orientation-preserving choice gives $f=g=0$, while the orientation-reversing coupling gives $f=g=2$.

Recall that our goal is to produce a coupling for which the particles are encouraged to couple in finite time.  Toward this end, it is not so much the specific coupling that matters to us as much as its qualitative features, especially the domination by a time-changed two-dimensional Bessel process.  In particular, the distance between the particles, which we called $r$, is a semi-martingale under any coupling, and this semi-martingale is characterized by the corresponding versions of $f$ and $g$.  We call any coupling with the feature that its corresponding versions of $f$ and $g$ satisfy all of the properties of $f$ and $g$ listed in the previous paragraph (that is, the domination by a time-changed two-dimensional Bessel process is strict expect for the configurations given, for which strictness fails in the ways described) an \emph{adequate} coupling.  As this terminology suggests, the coupling we've already described is the optimal one from our point of view, but any adequate coupling is still good enough for our exploration of minimal surfaces.  We will see why this matters in the next section.

\subsection{The martingale problem}\label{Section:Mart}

The computations in the previous section describe what the optimal coupling should be at a single point of $(x,y)\in M\times N\setminus\{r=0\}$; more specifically, they specify the desired covariance structure at each point.  The next task is to show that this gives a global description of the coupling.  (Dealing with existence is a separate question.  Here we just want to see that this pointwise choice of $\sigma$ and $A$ leads to a well-defined operator such that the corresponding martingale problem gives the coupling we're looking for.)  

As we have already indicated in the previous section, we think of a coupled Brownian motion as a process on the product space such that the $M$ and $N$ marginals are Brownian motions.  More specifically, we can formulate this in the language of the martingale problem, and we begin by briefly recalling what the martingale problem is (for a more detailed explanation, consult the standard text~\cite{SAndV}).  Let $\tilde{L}$ be a second-order operator on $M\times N$.  For simplicity, we consider the case when the solution is not allowed to explode (which will be the case for coupled Brownian motion on stochastically complete manifolds); extending this to allow for explosion is easily accomplished by stopping the process at the explosion time, as we did for Brownian motion at the beginning of Section~\ref{Section:Basic}.  Then using the notation from the beginning of Section~\ref{Section:Basic}, a solution to the martingale problem corresponding to $\tilde{L}$, starting at $(x_0,y_0)$, is a probability measure $P_{(x_0,y_0)}$ on $(C[0,\infty), \sB)$ such that 
\[\begin{split}
&P_{(x_0,y_0)}\lb \omega(0)=(x_0,y_0)\rb =1 \\
\text{and}\quad&\lp h(\omega(t))-\int_0^{t} \tilde{L} h(\omega(s))\, ds ; \sB_{t} \rp\quad \text{is a $P_{(x_0,y_0)}$-martingale,}
\end{split}\]
for any smooth, compactly supported function $h$ on $M\times N$.  Further, we note that the martingale problem is compatible with stopping times.  For example, if we let $\zeta$ be the first hitting time of $\{r=0\}$, then we can consider solutions to the martingale problem up until $\zeta$ by stopping the process at $\zeta$.  In this case, the operator $\tilde{L}$ need not be defined at $\{r=0\}$, and the class of test functions $h$ can be restricted to those having compact support on $(M\times N)\setminus\{r=0\}$.

Next, we review some basic facts about martingale problems that we will use throughout what follows.  There are two basic situations for which a solution to the martingale problem is known to exist.  If $\tilde{L}$ is $C^2$, then there is a unique solution to the corresponding martingale problem, at least until possible explosion (see Chapter 5 of~\cite{SAndV}).  If $\tilde{L}$ is locally bounded and uniformly elliptic, then there is a strong Markov family of solutions to the corresponding martingale problem, at least until possible explosion (see Exercises 7.3.2 and 12.4.3 of~\cite{SAndV}).  Further, both of these results can be localized (see Chapter 10 of~\cite{SAndV}).  Indeed, we have already taken advantage of this fact in transferring the above results to the manifold setting.  Localization also means that a global solution can be constructed by patching together local solutions.

As indicated, specifying a martingale problem corresponds to specifying a second-order operator, the diffusion operator.  At any point $(x,y)\in (M\times N)\setminus\{r=0\}$, we have product normal coordinates $(\alpha, \beta, a, b)$ as described in the previous section.  The operator at $(x,y)$ in the $(\alpha, \beta, a, b)$ coordinates is given by a symmetric, non-negative definite matrix $[a_{i,j}]|_{(x,y)}$ such that the operator applied to a function $f$ is given by $\frac{1}{2}\sum a_{i,j} \partial_i\partial_j f$.  Let
\[
\sO= \begin{bmatrix}
  A \cos\sigma & A\sin\sigma  \\ 
 -\sin\sigma & \cos\sigma  \end{bmatrix} ,
 \]
let $\sO^T$ be its transpose, and let $I$ be the $2\times 2$ identity matrix.
Then we have
\begin{equation}\label{Eqn:AatXY}\begin{split}
[a_{i,j}]|_{(x,y)}=&
 \begin{bmatrix} 1 & 0 & A \cos\sigma & -\sin\sigma \\ 
 0 & 1 & A \sin\sigma & \cos\sigma \\ 
 A \cos\sigma & A\sin\sigma & 1 & 0 \\ 
 -\sin\sigma & \cos\sigma & 0 & 1  \end{bmatrix} \\
 =& 
  \begin{bmatrix} 
 I & \sO^T  \\ 
 \sO & I \end{bmatrix}=
  \begin{bmatrix} I \\ 
 \sO  \end{bmatrix} \times
  \begin{bmatrix} I & \sO^T  \end{bmatrix} .
\end{split}\end{equation}
It's easy to see that $[a_{i,j}]|_{(x,y)}$ is symmetric and non-negative definite, and also that it has rank two.

Carrying out this procedure at every point shows that any pointwise choice of $A$ and $\sigma$ does in fact give rise to a well-defined operator.  Let $L_{+}$ be the operator determined by the optimal choice of coupling from the orientation-preserving component of $O(2)$, that is, the operator given by $(A,\sigma)=(1,\sigma_{+})$, and let $L_{-}$ be the analogous operator arising form the orientation-reversing component of $O(2)$, given by $(A,\sigma)=(-1,\sigma_{-})$.  Then define $\ol{L}$ to be equal to $L_{+}$ on the interior $\Sigma_{+}$ and $L_{-}$ on the interior of $\Sigma_{-}$.  This determines $\ol{L}$ uniquely except on $\Sigma_{0}$, where there are two possibilities.  We will have more to say about $\Sigma_{0}$ later, but for now we will leave open whether we choose the orientation-preserving or orientation-reversing possibility.  

Now suppose we have a solution to martingale problem corresponding to $\ol{L}$.  We want to compute the semi-martingale decomposition of any function $h$ on $M\times N\setminus\{r=0\}$ composed with the solution process.  The bounded variation process is just the integral of $Lh$ along paths, while the quadratic variation of the martingale part is given by the integral of $\nabla h \cdot [a_{i,j}]\cdot\nabla h$ along paths.  Alternatively, the growth of quadratic variation of the martingale part at each point (or more generally, the joint variation coming from two different functions) is determined by a bilinear form $\ol{\Gamma}(\cdot,\cdot)$.  In the same normal coordinates around $(x,y)$, we can write a vector $v\in T_{(x,y)}(M\times N)$ as 
\[
v=v_{\alpha}\partial_{\alpha}+v_{\beta}\partial_{\beta}+v_{a}\partial_{a}+v_{b}\partial_{b} .
\]
Then $\ol{\Gamma}(v,v)$ is determined by its expression in these coordinates.  This is most easily done by using the inner product structure on the tangent space.  Also, we divide the formula for $\ol{\Gamma}$ into two cases, the case where $(A,\sigma)=(1,\sigma_{+})$ corresponding to the diffusion operator at our point being $L_{+}$ and the case where $(A,\sigma)=(-1,\sigma_{-})$ corresponding to $L_{-}$.  We denote the corresponding bilinear forms by $\Gamma_{+}$ and $\Gamma_{-}$.  Then we have
\[\begin{split}
\Gamma_{+}(v,v) &= \ip{v}{ \partial_{\alpha} + \cos\sigma_{+} \partial_{a}-\sin\sigma_{+}\partial_{b}}^2 + \ip{v}{ \partial_{\beta} +\sin\sigma_{+}\partial_{a}+\cos\sigma_{+} \partial_{b}}^2 \\
\Gamma_{-}(v,v) &= \ip{v}{ \partial_{\alpha} -\cos\sigma_{-} \partial_{a}-\sin\sigma_{-}\partial_{b}}^2+ \ip{v}{ \partial_{\beta} -\sin\sigma_{-}\partial_{a}+\cos\sigma_{-} \partial_{b}}^2
\end{split}\]
Note that $\ol{\Gamma}$ at a point depends only on first-order information at that point.

We now verify that a solution to the martingale problem corresponding to $\ol{L}$ gives the coupling that it's supposed to.  First of all, any solution has marginals which are Brownian motions.  This follows from projecting $[a_{i,j}]|_{(x,y)}$ down to $M$ (which means taking the upper-right $2\times2$ sub-matrix), which shows that the marginal diffusion operator is expressed, in normal coordinates, as the $2\times 2$ identity matrix.  Thus, the marginal diffusion operator is one-half the Laplacian and the marginal process is Brownian motion.  The same holds for the marginal process on $N$.  Second, we can check that any solution induces a process $r_t=r(x_t,y_t)$ with the desired semi-martingale decomposition.  Again, choose $(\alpha,\beta,a,b)$ normal coordinates around a point $(x,y)$.  Since $r$ doesn't depend on a particular choice of Euclidean coordinates, we can choose $z_1$, $z_2$ and $z_3$ as at the beginning of the previous section.  The coordinate process will be martingales, and we see that $[a_{i,j}]|_{(x,y)}$ is chosen so that their covariance structure is as desired.  Once we know that the coordinate processes have the desired semi-martingale decompositions, the decomposition for $r_t$ follows.

Ideally, we would produce a solution to the martingale problem for $\ol{L}$, and use that in our study of minimal surfaces.  Unfortunately, $\ol{L}$ is everywhere degenerate and discontinuous at $\Sigma_0$.  Thus, the existence of a solution is not guaranteed by standard theorems.  More to the point, we don't have sufficient control of the configuration angles near $\Sigma_0$ to prove existence for arbitrary minimal surfaces.
To get around this difficulty, we will modify our operator.  The idea is that anywhere the coupling is strictly dominated by a time-changed two-dimensional Bessel process we have room to adjust it while remaining dominated by such a Bessel process.  In particular, this will allow us to produce an operator that is easier to deal with but such that any solution to the corresponding martingale problem will still be adequate (although no longer optimal, in the sense that a solution corresponding to $\ol{L}$ would be).

To make this precise, let $\Sigma_e\subset\Sigma_0$ be the set of points where $\psi=0$ (recall our convention that anywhere we can choose $\psi$ to be zero, we do) and $\theta+\phi=\pi/2$.  Thus $\Sigma_e$ is precisely the set of points in $\Sigma_0$ where the coupling looks instantaneously like a time-changed two-dimensional Bessel process, and where we have no room to perturb the operator while keeping the coupling adequate.  It follows that, for any $(x,y)\in \Sigma_0\setminus\Sigma_e$, the value of $f-g$ (recall that $f$ is the time-derivative of the quadratic variation and $g/2r$ is the time-derivative of the bounded variation part of $r_t$ under any solution to the martingale problem) achieved by the optimal coupling is positive.  Further, this value depends continuously on the cross-variations (which determine how the Brownian motions are coupled at a point), so a sufficiently small change in the operator will decrease $f-g$ (which we now view as depending on the operator as well) but maintain its positivity.  We can't make the operator continuous with an arbitrarily small perturbation, but we can make it elliptic.  In particular, we can replace the cross-variations prescribed by Equation~\eqref{Eqn:Coupling} (where we choose $\sigma$ and $A$ to correspond to the optimal coupling) with the cross-variations prescribed by
\[\begin{split}
d\CV{\alpha_t}{a_t} &= (1-\hat{\epsilon}) A\cos\sigma\, dt , \quad d\CV{\beta_t}{a_t}= (1-\hat{\epsilon}) A\sin\sigma\, dt, \\
d\CV{\alpha_t}{b_t} &= (1-\hat{\epsilon}) \lp-\sin\sigma\rp \, dt , \quad\text{and}\quad d\CV{\beta_t}{b_t}= (1-\hat{\epsilon}) \cos\sigma\, dt .
\end{split}\]
In terms of the corresponding operator characterizing the martingale problem, in $(\alpha,\beta,a,b)$ normal coordinates at a point, it is given by the following matrix,
\[
\begin{bmatrix} I & \sqrt{1-\hat{\epsilon}^2} \sO^T \\
\sqrt{1-\hat{\epsilon}^2} \sO & I 
\end{bmatrix}=
\begin{bmatrix} I & 0 \\
\sqrt{1-\hat{\epsilon}^2} \sO & \hat{\epsilon} I 
\end{bmatrix} \times
\begin{bmatrix} I & \sqrt{1-\hat{\epsilon}^2} \sO^T \\
0 & \hat{\epsilon} I 
\end{bmatrix} ,
\]
where $I$ and $\sO$ are as in Equation~\eqref{Eqn:AatXY}.  Heuristically, from the stochastic differential equation point of view, this means that Brownian motion on $N$ is only partially being driven by Brownian motion on $M$, and ``the rest of'' the Brownian motion on $N$ is being driven by an independent source, where the ratio between these two driving sources is governed by $\hat{\epsilon}$.  In particular, $\hat{\epsilon}=0$ gives us the same coupling as before, while $\hat{\epsilon}=1$ means that the Brownian motions on $M$ and $N$ are independent.  For any $\hat{\epsilon}\in[0,1]$, the resulting marginals are still Brownian motions on $M$ and $N$.  Also, the value of $f-g$ achieved by the coupling is smooth in $\hat{\epsilon}$.   Finally, note that if $\hat{\epsilon}<1$, then the corresponding operator will be elliptic.

It follows that we can choose $\hat{\epsilon}\in[0,1/2)$ to depend smoothly on $(x,y)\in M\times N\setminus\{r=0\}$ in such away that $\hat{\epsilon}=\hat{\epsilon}(x,y)$ is zero outside of some neighborhood of $\Sigma_0$, zero on $\Sigma_e$, positive on $\Sigma_0\setminus\Sigma_e$, and so that the corresponding coupling is adequate (which also requires $\hat{\epsilon}$ to be zero on the set of points where $\psi=0$ and $\theta=\phi<\pi/4$).  Call the resulting operator $L$.  Of course, $L$ depends on our precise choice of a function $\hat{\epsilon}$.  We will simply assume that we've chosen some $\hat{\epsilon}$ satisfying the criteria just mentioned, and any will do.  It follows from our construction that, for any point $(x,y)\in \Sigma_0\setminus\Sigma_e$, $L$ is uniformly elliptic in some neighborhood of $(x,y)$, and also that $L$ has the same smoothness properties as $\ol{L}$ everywhere.  The advantage of using $L$ instead of $\ol{L}$ is that now we really only need to worry about the behavior of the process near the smaller set of discontinuities $\Sigma_e$, which we will find much more manageable.

It is the martingale problem for $L$ for which we will show there exists a solution.  As mentioned, this will be sufficient for our purposes, since the resulting coupling will be adequate.  Finally, we comment that our study of $\ol{L}$ remains relevant, since we will primarily approach $L$ as a perturbed version of $\ol{L}$.  This is especially true at $\Sigma_e$, since the two operators are equal there and thus we will continue to use our expressions for $\ol{\Gamma}$.

\section{Existence of an adequate coupling}

\subsection{Preliminary results}

The fact that we have defined $\ol{L}$, and hence $L$, in a different system of coordinates at each point obscures the questions of how smooth $L$ is.  This isn't too hard to get around near points where the configuration $(\theta,\phi,\psi)$ is in the interior of its range, since one can see that all of the coordinates (both on $\Re^3$ and on $M\times N$) introduced in the previous section will be smooth near such a point.  However, these canonical coordinate systems can be discontinuous when the configuration is at the boundary if its range.  This means, for example, that the identification of $O(2)$ with isometries between the tangent spaces can also change discontinuously.

One could try to introduce other, better behaved coordinates near such points.  However, this turns out not to be necessary.  The following lemma clarifies the behavior of $L$.

\begin{Lemma}\label{Lemma:SmoothL}
Using the above notation, $\Sigma_0$ is locally the zero level-set of a smooth function, and $L$ is smooth at any point not in $\Sigma_0$.
\end{Lemma}

\emph{Proof:}  As already mentioned, $L$ is smooth at a point if and only if $\ol{L}$ is.  Thus we will prove the lemma for $\ol{L}$, since we have more explicit formulas for $\ol{L}$.

We begin by considering the expression of an operator in local, smooth coordinates.  In particular, choose smooth product coordinates $(x,y)$ on some product neighborhood $S=S_M\times S_N$ of a point $(x_0,y_0)\in (M\times N)\setminus\{r=0\}$, and smooth orthonormal frames for both $M$ and $N$ on $S_M$ and $S_N$.  As above, this choice of orthonormal frames on $M$ and $N$ determines an identification of the isometries between tangent planes and $O(2)$.  Suppose the operator under consideration is determined by a choice of element of $O(2)$ at each point $(x,y)\in S$, in the fashion described above.  Then if the the map from $S$ to $O(2)$ is smooth, so is the corresponding diffusion operator.  This is just a consequence of the change of variables formula and the fact that all of the functions involved are smooth.

Thus, in order to show that $\ol{L}$ is smooth at such a point $(x_0,y_0)$, it is enough to show that, with respect to coordinates and frames as above, the map from $S$ into $O(2)$ which determines $\ol{L}$ is smooth.  Suppose that, with respect to the induced identification of the isometries between tangent planes and $O(2)$, $\ol{L}$ is determined by elements of the orientation-preserving component of $O(2)$ on $S$.  Let $s\in[0,2\pi)$ be the natural coordinate on the orientation-preserving component.  Then Equation~\eqref{Eqn:FMinusG} implies that
\[
(f-g)(x,y,s) = \alpha(x,y)\cdot\cos s +\beta(x,y)\cdot \sin s +\gamma(x,y)
\]
for some functions $\alpha$, $\beta$, and $\gamma$.  Here $f$ and $g$ are the functions determining the semi-martingale decomposition of $r_t$, under the coupling corresponding to any choice of $s$ (assuming such a coupling exists).  In particular, this equation gives $f-g$ for any choice of $s$; the value of $f-g$ for $\ol{L}$ is obtained by letting $s$ be the ``optimal choice'', as determined above, at each $(x,y)$.

Note that the images of $x$ and $y$ in $\Re^3$ and the corresponding tangent planes vary smoothly, as do the projections of the tangent planes onto the direction of $x-y$ (since $r\neq 0$) and its orthogonal complement.  Since $f$ and $g$ depend smoothly on these projections and on the coupling, it follows that $f-g$ is a smooth function of $(x,y,s)$.  Further, Equation~\eqref{Eqn:FMinusG} makes it clear that $\gamma(x,y)$ is always non-positive, and zero only when $(x,y)$ is such that $\theta=\phi=0$ in the standard configuration.  Since we showed in the previous section that the maximum of $f-g$ is always non-negative, and since this maximum is equal to $\sqrt{\alpha^2+\beta^2}+\gamma$ on $S$ (recall that we're assuming that the maximum is achieved by a coupling in the orientation-preserving component on $S$), it follows that $\sqrt{\alpha^2+\beta^2}$ is positive except possibly where $\theta=\phi=0$.  Where $\theta=\phi=0$, this can be checked by hand, and thus $\sqrt{\alpha^2+\beta^2}$ is positive everywhere in $S$.

Next, observe that because $f-g$ is smooth in $x$ and $y$ for every fixed $s\in[0,2\pi)$, we can see that $\alpha$, $\beta$, and $\gamma$ are smooth functions of $x$ and $y$.  This implies that $\sqrt{\alpha^2+\beta^2}$ is positive and smooth on $S$, and thus, possibly by shrinking $S$, we can assume that $\sqrt{\alpha^2+\beta^2}$ is smooth and uniformly positive on $S$.  We know that, for given $x$ and $y$, the choice of $s$ that maximizes $f-g$, which we will denote $s_+(x,y)$, is the angular component of the polar representation of $(\alpha(x,y),\beta(x,y))\in\Re^2$.  Since the components of this vector are smooth and and the length of the vector is uniformly positive on $S$, it follows that $s_+$ is smooth on $S$.  Since $s_+$ determines $\ol{L}$ on $S$, it follows that $\ol{L}$ is smooth on $S$.  Obviously, an analogous argument applies if we assume that $\ol{L}$ is realized by an element of the orientation reversing component of $O(2)$ on $S$.

Now suppose that we choose a point $(x_0,y_0)\in\Sigma_0$, and we again choose smooth coordinates and smooth orthonormal frames on a neighborhood $S$ as above.  At $(x_0,y_0)$, there are two possible choices of optimal couplings, one from each connected component of $O(2)$.  Let $\ol{L_{+}}$ be the operator on $S$ obtained from the optimal choice of coupling in the orientation-preserving component, and let $\ol{L_{-}}$ be the analogous operator for the orientation-reversing component (these need not the be the same as $L_{+}$ and $L_{-}$, as these were defined for a particular choice of orthonormal frames).  Since the value of $f-g$ achieved by each of these two couplings is positive at $(x_0,y_0)$, the above argument implies that, after possibly shrinking $S$, we can assume that both $\ol{L_{+}}$ and $\ol{L_{-}}$ are smooth on $S$.  Further, the values of $f-g$ realized by $\ol{L_{+}}$ and $\ol{L_{-}}$ are also smooth on $S$.  Since $\Sigma_0$ is the set where these two values are equal, it follows that $\Sigma_0$ is the zero-level set of a smooth function on $S$.

Because $\Sigma_0$ is locally the zero level-set of a smooth function, it is a closed set.  Thus, if we choose a point $(x_0,y_0)$ in the complement of $\Sigma_0$, we can take our neighborhood $S$ to be disjoint from $\Sigma_0$.  Then the optimal coupling belongs to the same connected component of $O(2)$ throughout $S$ (that the coupling can only ``switch'' components of $O(2)$, with respect to smooth frames, at $\Sigma_0$ follows from the smoothness of $f-g$ and the fact that $\Sigma_0$ is the only set where there are two possible choices for the optimal coupling), and the above argument shows that $\ol{L}$ is smooth on $S$.  $\Box$

This proof also gives a description of the discontinuities of $\ol{L}$ and $L$.  We know that $\ol{L}$ can only be discontinuous at $\Sigma_0$.  In a neighborhood of a point in $\Sigma_0$, both $\ol{L_{+}}$ and $\ol{L_{-}}$ (using the notation from the proof, and assuming that we have chosen some smooth orthonormal frames on $M$ and $N$) are smooth.  The discontinuity in $\ol{L}$ arises from the fact that we ``switch'' from having $\ol{L}$ given by $\ol{L_{+}}$ to having $\ol{L}$ given by  $\ol{L_{-}}$.  In particular, if we let $h$ be a smooth function such that $\Sigma_0$ is locally the zero level-set of $h$, and if the gradient of $h$ is non-zero at $(x_0,y_0)\in\Sigma_0$, then $\Sigma_0$ is locally a hypersurface, and $\ol{L}$ is given by $\ol{L_{+}}$ or $\ol{L_{-}}$ depending on which side of the hypersurface we're on, at least locally.  An analogous description applies to $L$, with $\ol{L_{+}}$ and $\ol{L_{-}}$ replaced by their perturbed versions.

One consequence of this lemma is that if we start our coupled Brownian motion, corresponding to $L$, in the complement of $\Sigma_{0}$ smoothness implies that we have a unique solution at least until the first time the process hits $\Sigma_{0}$.  Further, the ellipticity of $L$ gives us what we need for existence on $\Sigma_0\setminus\Sigma_e$, which is the reason for introducing $L$.  Near $\Sigma_e$, we need the following lemma.  Recall that the signed distance to a smooth hypersurface is the smooth function, defined in some neighborhood of the hypersurface, the absolute value of which is the distance to the hypersurface.

\begin{Lemma}\label{Lemma:SigmaE}
Let $(x,y)$ be a point of $\Sigma_e$.  Then there is a neighborhood $S$ of $(x,y)$ such that one of the following holds:
\begin{enumerate}
\item $\Sigma_0 \cap S$ is a smooth hypersurface, and if $v_0$ is the gradient of the signed distance to $\Sigma_0$, then either $\Gamma_{+}(v_0,v_0)$ or $\Gamma_{-}(v_0,v_0)$ is positive at $(x,y)$.
\item $\Sigma_e \cap S$ is contained in some smooth hypersurface $H$, and if $v_H$ is the gradient of the signed distance to $H$, then both $\Gamma_{+}(v_H,v_H)$ and $\Gamma_{-}(v_H,v_H)$ are positive at $(x,y)$.
\end{enumerate}
\end{Lemma}

\emph{Proof:}  We begin by introducing some notation.  We let
\[
\ol{\partial_{\alpha}} = \partial_{\alpha} + A\cos\sigma \partial_{a}-\sin\sigma\partial_{b}
\quad\text{and}\quad \ol{\partial_{\beta}} = \partial_{\beta} +A\sin\sigma\partial_{a}+\cos\sigma \partial_{b} ,
\]
so that $\Gamma_{\pm}(v,v) = \ip{v}{\ol{\partial_{\alpha}}}^2 + \ip{v}{\ol{\partial_{\beta}}}^2$ for the appropriate choice of $A$ and $\sigma$.

We need to to look at the first order derivatives of $\theta$, $\phi$, and $\psi$ with respect to $\ol{\partial_{\alpha}}$ and $\ol{\partial_{\beta}}$ when $\psi=0$ and $\theta+\phi=\pi/2$.  Fortunately, the derivatives of $(x-y)/|x-y|$ are particularly simple in this case, as we see from Equation~\eqref{Eqn:SystemEvo} and the fact that both $\sigma_{+}$ and $\sigma_{-}$ are zero when $\phi=0$.  As for the normal vectors $m$ and $n$, their derivatives are constrained only by the fact that the Gauss map is anti-conformal.  Thus, we can describe the derivatives of $m$ by $k_1\geq 0$ and $s_1\in[0,2\pi)$ in the following way.  We know that $\partial_{\alpha}$ and $\partial_{\beta}$ form an orthonormal basis for $T_xM$ and thus also for $T_{m}\S^2$.  Then the general anti-conformal map between the two can be written as
\[
\partial_{\alpha} (m) = k_1\lp \cos s_1\partial_{\alpha} + \sin s_1\partial_{\beta}\rp\quad \text{and} \quad \partial_{\beta}(m) = k_1 \lp \sin s_1\partial_{\alpha} -\cos s_1 \partial_{\beta} \rp .
\]
An analogous description of the derivatives of $n$ can be given in terms of $k_2\geq 0$ and $s_2\in[0,2\pi)$ with respect to $\partial_a$ and $\partial_b$.

Because the roles of $M$ and $N$ are symmetric, we can assume that $\phi\geq \theta$.  Recall that $\Sigma_0$ is the zero level-set of $\cos\theta\cos\phi - \cos\psi\sin\theta\sin\phi$, at least when the configuration is in the interior of its range.  Next, note that this continues to hold if we allow $\psi$ to take negative values.  Thus, if we don't insist on using only configurations in the canonical range mentioned above, all three angles are smooth functions near $(x,y)$, if we assume that $\theta>0$, and $\Sigma_0$ is the zero level-set of a smooth function.  So we first consider the case where $\theta>0$.  In light of the above discussion, a few simple computations show that under these conditions, that is whenever $\psi=0$, $0<\theta\leq\pi/4$, and $\phi=\pi/2-\theta$, we have
\[\begin{split}
&\ol{\partial_{\alpha}}\lp\theta+\phi\rp = \frac{2}{r}\lp \cos\theta -A\sin\theta\rp-k_1\cos s_1-Ak_2\cos s_2 ; \\
&\ol{\partial_{\beta}}\lp\theta+\phi\rp = -k_1\sin s_1 -k_2\sin s_2 ; \\
&\ol{\partial_{\alpha}}\psi= -\frac{1}{\theta}k_1\sin s_1 + \frac{A}{\frac{\pi}{2}-\theta}k_2\sin s_2 ; \quad\text{and}\quad \ol{\partial_{\beta}}\psi= \frac{1}{\theta}k_1\cos s_1 - \frac{1}{\frac{\pi}{2}-\theta}k_2\cos s_2 .
\end{split}\]

We begin by determining when the first possibility in the lemma holds.  It's easy to see that $v_0$ at $(x,y)$ will be a non-zero multiple of the gradient of $\theta+\phi$.  As long as $\ol{\partial_{\alpha}}\lp\theta+\phi\rp$ for at least one choice of $A$ or $\ol{\partial_{\beta}}\lp\theta+\phi\rp$ is not zero, either $\Gamma_{+}(v_0,v_0)$ or $\Gamma_{-}(v_0,v_0)$ will be positive at $(x,y)$.  We conclude that the first condition holds unless all three of the following equations are satisfied (at $(x,y)$):
\[
-k_1\sin s_1 =k_2 \sin s_2 , \quad
\frac{2}{r}\cos\theta = k_1\cos s_1 ,\quad\text{and }
\frac{2}{r}\sin\theta = -k_2\cos s_2 .
\]
To complete the proof, for the $\theta>0$ case, we need to show that the second possibility holds whenever all three of the equations are satisfied.  We know that $\Sigma_e$ is contained in the zero level-set of $\psi$, at least near $(x,y)$.  Note that, if the above three equations are satisfied, then $k_1\cos s_1>0$ and $k_2\cos s_2<0$.  This implies that $\ol{\partial_{\beta}}\psi>0$, for both choices of $A$, since it doesn't depend on $A$.  Thus we can let $H$ be the zero level-set of $\psi$, and the second possibility holds.

We now consider a point $(x,y)\in\Sigma_e$ where $\theta=0$.  Here, $\psi$ is not continuous, and $\theta$ (which, we recall, is the distance in $\S^2$ between $m$ and $(x-y)/|x-y|$) has the usual non-differentiability of a distance function at its zero level-set.
First, observe that $\psi=0$ requires that $m$, $n$, and $(x-y)/|x-y|$ all lie on the same geodesic in $\S^2$.  If we let $\gamma$ be the great circle through $(x-y)/|x-y|$ and $n$, then $\gamma$ varies smoothly near $(x,y)$.  Further, if we let $h$ be the signed distance between $m$ and $\gamma$, then $h$ is smooth near $(x,y)$ as well.  It follows that $\Sigma_e$ is contained in the zero level-set of $h$, near $(x,y)$.  Because of this, it is natural to ask what conditions must be satisfied in order for the second possibility in the lemma to hold.  To first order, $(x-y)/|x-y|$ will only move along $\gamma$.  It follows that if $k_1\neq 0$, either $\ol{\partial_{\alpha}}h$, for both values of $A$, or $\ol{\partial_{\beta}}h$ will be non-zero.  In this case, possibility two of the lemma holds, with $H$ the zero level-set of $h$.

If $k_1=0$, then $h$ is zero to first order.  In other words, $m$, $n$, and $(x-y)/|x-y|$ remain colinear (in the sense of lying on the same great circle) to first order.  This suggests the following approach.  If we allow $\phi$ to take values greater than $\pi/2$ and $\theta$ to take negative values (this is well-defined if we restrict $\theta$ to the great circle), then these two angles determine the configuration, up to first order.  It's easy to check that the characterization of $\Sigma_0$ as points with $\theta+\phi=\pi/2$ (when $h=0$) extends to these values of $\phi$ and $\theta$.  Further, both angles are differentiable and our earlier formulas specialize to
\[
\ol{\partial_{\alpha}}\lp\theta+\phi\rp = \frac{2}{r} -Ak_2\cos s_2 \quad\text{and}\quad 
\ol{\partial_{\beta}}\lp\theta+\phi\rp = -k_2\sin s_2 .
\]
It is immediate that these cannot both be zero for both values of $A$.  Since we know a priori, by Lemma~\ref{Lemma:SmoothL}, that $\Sigma_0$ is the zero level-set of a smooth function, it follows that the first possibility in the lemma holds.  $\Box$

The next two lemmas concern the case of an operator on $\Re^d$ which, in standard coordinates $(z_1,\ldots,z_d)$, can be written in terms of a (measurable) locally bounded function $a$ taking values in the set of symmetric, non-negative definite matrices and a (measurable) locally bounded function $b$ taking values in $\Re^d$ as
\be\label{StandardOp}
\tilde{L} = \frac{1}{2}\sum_{i,j=1}^d a_{i,j}(z) \partial_{z_i} \partial_{z_j} +\sum_{i=1}^d b_i(z) \partial_{z_i} ,
\ee
where the $a_{i,j}$ and the $b_i$ are the obvious components of $a$ and $b$.  Further, we choose $C>0$ so that $|a_{i,j}|$ and $|b_{i}|$ are all less than $C$ on $B_R= \{z: z_1+\cdots+z_n\leq R\}$.  Let $\tilde{P}$ be a measurable, strong Markov family of solutions to the martingale problem associated to $\tilde{L}$, with $\tilde{P}_{x}$ denoting the solution started at $x$ (we assume that such a family of solutions exists).

\begin{Lemma}\label{Lemma:Elliptic}
Let $\tilde{L}$ and $\tilde{P}_x$ be as above, and suppose further that $a$ is uniformly elliptic on $B_R$ (so that all of its eigenvalues are bounded below by some $c>0$).  Then the expected occupation time under $\tilde{P}_x$ from time 0 to $T$ has a density $G_T(x,y)$ on $B_{R/2}$, and this density obeys the estimate
\[
\sup_{x\in \Re^d} \|G_T(x,\cdot)\|_{L^q(B_{R/2})} \leq A
\]
where $A$ and $q$ are positive numbers, with $q>d/(d-1)$, that depend only on
$d$, $R$, $T$, $c$, and $C$.
\end{Lemma}
\emph{Proof:} Because $a$ and $b$ are locally bounded, a standard localization argument implies that it is enough to prove the lemma when the process is stopped at the first exit time of $B_R$.

Let $\sigma_1$ be the first exit time of $B_R$ (which is zero if $x\not\in B_R$); we denote the corresponding density by $G_{T\wedge\sigma_1}(x,y)$.  An application of Girsanov's theorem (see section 6.4 of~\cite{SAndV}) shows that it's enough to consider the case when $b\equiv 0$.  Then the estimate for $G_{T\wedge\sigma_1}(x,y)$, as well as the fact that such a density exists, follows directly from Corollary 2.4 of~\cite{FabesStroock}.  $\Box$

For the next lemma, we assume that the rectangle $(0,1)\times(-1,1)^{d-1}$ is contained in $B_R$.  Because we can rescale the coordinates, this is no loss of generality.

\begin{Lemma}\label{Lemma:Occupation}
Let $\tilde{L}$ and $\tilde{P}_x$ be as above.  Suppose that $a_{1,1}$ is bounded from below by a positive constant $c$ on the rectangle $(0,1)\times(-1,1)^{d-1}$.  Then the expected occupation time of $(0,\delta)\times(-1/2,1/2)^{d-1}$ from time 0 to $T$ under $\tilde{P}_x$ goes to zero with $\delta$, at a rate which depends only on $d$, $R$, $T$, $c$, and $C$.
\end{Lemma}

\emph{Proof:}  Again, the local boundedness of $a$ and $b$ means that it is sufficient to prove the lemma when the process is stopped at $\sigma_1$, the first exit from $B_R$.

Let $\xi(x)$ be a smooth, non-negative, even function such that $\xi(x)\leq |x|$, $|\xi'(x)|\leq 1$, $\xi''(x)\geq 0$, and $\xi''(x)=1$ on $(-1/2,1/2)$.  Then let $\xi_{\delta}(x) = \delta^2\xi(x/\delta)$.  Because $\xi_{\delta}(z_1-\delta/2)$ is bounded on $B_R$, its expectation at time $T\wedge\sigma_1$ is bounded from above by $k\delta$ for some constant $k$ depending only on $d$ and $R$ (assuming $\delta<1$).  On the other hand, It\^o's rule implies that the expectation of $\xi_{\delta}(z_1-\delta/2)$ at time $T\wedge\sigma_1$ is at least $c/2$ times the expected occupation time of $(0,\delta)\times(-1/2,1/2)^{d-1}$, minus $(T\wedge\sigma_1)C\delta$.  It follows that the expected occupation time of $(0,\delta)\times(-1/2,1/2)^{d-1}$ is less than
\[
\frac{ k+C(T\wedge\sigma_1) }{c/2}\delta .
\]
Thus we've proved that the lemma holds for the process stopped at $\sigma_1$.  $\Box$

\subsection{Proof of existence}

We are now in a position to prove the existence of an adequate coupling, namely, the existence of a solution to the martingale problem corresponding to any choice of operator $L$ as described above.

\begin{THM}\label{THM:Main}
Let $M$ and $N$ be any stochastically complete minimal surfaces.  For any points $x_0\in M$ and $y_0\in N$ with $r(x_0,y_0)\neq0$, there exists an adequate coupling of Brownian motions started at $x_0$ and $y_0$ defined until the first time $r_t=r(x_t,y_t)$ hits zero, in particular, the coupling corresponding to any choice of operator $L$ as described above.  Further, given any such $L$, the corresponding coupling is unique until the first time it hits $\Sigma_0$.
\end{THM}

\emph{Proof:}  Existence of a solution to the martingale problem for $L$ starting from $(x_0,y_0)$ is not guaranteed by standard results, so we proceed by an approximation argument.  Consider the family of operators $L^{(j)}=L +\Lap/j$ where $\Lap$ is the Laplacian on the product manifold $M\times N$.  In order to deal with stopping the process at $\{r=0\}$, we will need a second level of approximation.  For all small, positive $\epsilon$, let $\eta_{\epsilon}$ be a smooth function from $M\times N$ into $[0,1]$ that is 1 on $\{r\leq \epsilon/2\}$ and 0 on $\{r\geq \epsilon \}$.  Then we let
\[
L^{(j,k)}=(1-\eta_{1/k}) L^{(j)} + \eta_{1/k} \lp1+\frac{1}{j}\rp\Lap .
\]
We note that $L^{(j,k)}$ is defined on all of $M\times N$, unlike $L$ which is not defined on $\{r=0\}$.  

In order to apply standard theorems in martingale theory, it will be helpful to introduce global coordinates and thus work on $\Re^4$.  As usual, we can pass to the universal covers of $M$ and $N$, and so we assume that they are simply connected.  Because they have non-positive curvature, normal coordinates around any point give global coordinates, and map each surface diffeomorphically to $\Re^2$.  This extends to the product manifold in the obvious way.  In these global coordinates, each $L^{(j,k)}$ can be written in the form shown in Equation~\eqref{StandardOp} with locally bounded, locally uniformly elliptic coefficients.  It follows that for each $L^{(j,k)}$ we can find a measurable, strong Markov family of solutions to the corresponding martingale problem (see the beginning of Section~\ref{Section:Mart} for the definition of the martingale problem) which we denote $P^{(j,k)}_{(x_0,y_0)}$.  (A priori, these solutions are defined only up to explosion, but we will see in a moment that they never explode).  In order to simplify the notation, we will assume some starting point $(x_0,y_0)$ has been chosen and simply denote the corresponding measures by $P^{(j,k)}$ whenever there is no possibility of confusion.  Consider the marginal processes on $M$ and $N$.  In each case, the marginal distribution of $P^{(j,k)}$ solves the martingale problem for $(j+2)/2j$ times the Laplacian.  This is just time-changed (by a constant factor) Brownian motion on a stochastically complete manifold.  It follows that the $P^{(j,k)}$ processes never explode, since $M$ and $N$ are both stochastically complete and a process on a product manifold blows up if and only if one of the marginals blows up.

Next, we need to show that any sequence of $\{P^{(j,k)}\}$ has a weakly convergent subsequence, that is, that this family of measures is pre-compact.  All of our process start from the same point, so Theorem 1.3.1 of~\cite{SAndV} asserts that $\{P^{(j,k)}\}$ is pre-compact if and only if, for every $\rho>0$ and $T<\infty$, we have
\begin{equation}\label{Eqn:PreCompact}
\lim_{\delta\searrow 0} \inf_{j,k} P^{(j,k)}
\lb \sup_{\substack{0\leq s\leq t\leq T \\ t-s\leq\delta}} \lab \omega(t) - \omega(s) \rab_{\Re^4}\leq\rho\rb=1 .
\end{equation}
Recall that the marginals on $M$ and $N$ are time-changed (by a constant factor) Brownian motions, and thus we see that the family of marginals on $M$ and the family of marginals on $N$ both possess the property described by Equation~\eqref{Eqn:PreCompact}.  For any path on the product space, the increment $\lab \omega(t) - \omega(s) \rab_{\Re^4}$ is bounded by the sum of the increments of the projections onto $M$ and $N$, by the triangle inequality.  It follows that the family $\{P^{(j,k)}\}$ possesses the property described by Equation~\eqref{Eqn:PreCompact}, and thus any sequence of $\{P^{(j,k)}\}$ has a weakly convergent subsequence.

Consider any sequence $(j(l),k(l))$ such that $j(l)\ra\infty$ and $k(l)\geq k_0$ for some positive integer $k_0$.  The corresponding sequence of measures has a convergent subsequence, so after possibly re-indexing our sequence, we can assert that $P^{(j(l),k(l))}$ converges to a limit we call $P^{k_0}$  (obviously, the limit depends in general on the sequence $(j(l),k(l))$ and not just on $k_0$, but we'll see that this notation will be sufficient for our purposes).  Further, let $\zeta_{\epsilon}$ be the first hitting time of the set $\{r\leq \epsilon\}$; in particular, $\zeta_0=\zeta$ which we have previously defined as the first hitting time of $\{r=0\}$.

We wish to prove that $P^{k_0}$ is a solution to the martingale problem corresponding to $L$ until $\zeta_{1/k_0}$.  It's easy to see that $P^{k_0} \lb \omega(0)=(x_0,y_0)\rb =1$ almost surely, and so it remains to prove that $P^{k_0}$ has the desired martingale property.  For this, it is enough to show that
\[
\E\lb F \lp  h(\omega(t))-h(\omega(s))-\int_{s}^{t} L h(\omega(u))\, du \rp\rb =0
\]
for $0\leq s< t$, any bounded, continuous, $\sB_{s}$-measurable function $F:C[0,\infty)\ra \Re$, and any smooth $h$ compactly supported on $(M\times N)\setminus\{r\leq 1/k_0\}$, where the expectation is with respect to $P^{k_0}$.  This in turn will follow if we show that
\be\label{Eqn:Martingale}\begin{split}
& \E^{(l)}\lb |F|  \int_{s}^{t} \lab \lp L^{(j(l))}-L\rp h(\omega(u))\rab \, du \rb \\
&\quad\quad\quad +  \lab \E^{(l)}\lb F\int_{s}^{t} Lh(\omega(u))\, du\rb - \E\lb F\int_{s}^{t} Lh(\omega(u))\, du\rb\rab 
\end{split}\ee
goes to zero as $l$ goes to infinity, where $\E^{(l)}$ is expectation with respect to $P^{(j(l),k(l))}$.  The first term goes to zero because $F$ is bounded and $\lp L^{(j(l))}-L\rp h$ converges to zero uniformly, using the fact that $h$ is smooth and compactly supported on $(M\times N)\setminus\{r\leq 1/k_0\}$ and $L^{(j(l))}$ converges to $L$ on $(M\times N)\setminus\{r\leq 1/k_0\}$.

For the second term, we note that, by a partition of unity argument, it is sufficient to show that for any point $(x,y)\in\lp M\times N\rp\setminus\{r\leq 1/k_0\}$ there is an open neighborhood $S$ of the point such that the second term goes to zero for $h$ supported on $S$ (and thus $S$ should be taken to be disjoint from $\{r\leq 1/k_0\}$).  For the purposes of such an argument, we see that there are three types of points.  First, suppose that $(x,y)$ is contained in the complement of $\Sigma_0$.  Then we can choose $S$ also to be contained in this interior.  If $h$ is smooth and supported on $S$, then so is $Lh$, using the fact that $L$ is smooth on the complement of $\Sigma_0$.  Then the second term of Equation~\eqref{Eqn:Martingale} goes to zero by the definition of weak convergence (recall that $F$ is continuous).

Next, suppose that $(x,y)\in \Sigma_0\setminus\Sigma_e$.  Then $L$ is uniformly elliptic and bounded on any sufficiently small neighborhood of $(x,y)$.  Moreover, we can choose coordinates centered at $(x,y)$ so that $L^{(j(l))}$ and $P^{(j(l),k(l))}$ satisfy the assumptions of Lemma~\ref{Lemma:Elliptic}  for some choice of constants $R$, $c$, and $C$ independent of $n$, and with $d=4$ and $T>t$.  Also, we can choose $S$, our neighborhood of $(x,y)$, to be contained in $B_{R/2}$ (in these coordinates) and disjoint from $\{r\leq 1/k_0\}$.  Let $q>4/3$ be the constant from Lemma~\ref{Lemma:Elliptic}, and let $q'$ be its H\"{o}lder conjugate.  Then we can find a sequence of continuous functions $\xi_{m}$ supported on $B_{R/2}$ which approximate $L h$ in $L^{q'}(B_{R/2})$.  By weak convergence, we know that
\[
\E^{(l)}\lb F \int_{s}^{t} \xi_m(\omega(u))\, du \rb \ra
\E\lb F \int_{s}^{t} \xi_m(\omega(u))\, du \rb 
\]
as $l \ra\infty$.  Further, Lemma~\ref{Lemma:Elliptic} and the H\"older inequality imply that
\[
\E^{(l)}\lb F \int_{s}^{t} \xi_m(\omega(u))\, du \rb \ra
 \E^{(l)}\lb F \int_{s}^{t} Lh(\omega(u))\, du \rb
\]
as $m\ra\infty$, uniformly in $l$.  Because this convergence is uniform in $l$, combining these two equations shows that the second term of Equation~\eqref{Eqn:Martingale} goes to zero as desired.

The final type of point we need to consider is $(x,y)\in\Sigma_e$.  This divides further into two cases, depending on which of the possibilities in Lemma~\ref{Lemma:SigmaE} holds.  Suppose that the second possibility in Lemma~\ref{Lemma:SigmaE} holds.  
Let $\Gamma$ be the bilinear form that gives the cross-variation of vector fields under $L$.  We know that the operators corresponding to the optimal orientation-preserving and reversing couplings are smooth, and that $L$ agrees with $\ol{L}$ at $(x,y)$.  It follows that $\Gamma(v_H,v_H)$ is bounded below by a positive constant on an open neighborhood $S'$ of $(x,y)$; we also assume that $S'$ is disjoint from $\{r\leq 1/k_0\}$.  
Thus, we can find coordinates $(z_1,\ldots,z_4)$ centered at $(x,y)$ such that $L^{(j(l))}$ and $P^{(j(l),k(l))}$ satisfy the assumptions of Lemma~\ref{Lemma:Occupation} on both sides of $H$ for some choice of constants $R$, $c$, and $C$ independent of $l$, with $d=4$ and $T>t$; and where these coordinates are such that 
\[ 
(-1,1)^4\subset S' \subset B_R 
\]
and $z_1$ restricted to $S'$ is the signed distance to $H$.  In particular, Lemma~\ref{Lemma:Occupation} implies that the occupation time of $(-\delta,\delta)\times(-1/2,1/2)$ goes to zero in $\delta$ at a rate that can be taken to be independent of $l$.

We can assume that the support of $h$ is contained in $S= (-1/2,1/2)^4$. For any $\delta\in(0,1/2)$, we can find a 
mollified version of $Lh$, which we denote $\xi_{\delta}$, with the properties that $\xi_{\delta}$ is supported on $S$, $\xi_{\delta}$ is continuous in a $\delta/2$-neighborhood of $\Sigma_e$, the $\xi_{\delta}$ are bounded uniformly in $\delta$, $Lh$ and $\xi_{\delta}$ are equal outside of a $\delta$-neighborhood of $\Sigma_e$, and the discontinuities of $\xi_{\delta}$ are contained in an open set where $L$ is uniformly elliptic (here we mean uniformly in space, not in $\delta$).  We claim that, for any $\delta$,
\be\label{Eqn:MoreWeakCon}
\E^{(l)}\lb F \int_{s}^{t} \xi_{\delta}(\omega(u))\, du \rb \ra
\E\lb F \int_{s}^{t} \xi_{\delta}(\omega(u))\, du \rb 
\ee
as $l\ra\infty$.  To see this, note that a bump function argument shows that $\xi_{\delta}$ can be written as the sum of a bounded function supported on a subset of $S$ where $L$ is uniformly elliptic and a continuous function supported on $S$.  The argument given above for $Lh$ supported on a set where $L$ is uniformly elliptic shows that we have the desired convergence for the first term in the decomposition of $\xi_{\delta}$, and weak convergence applies directly to the second term in the decomposition.  The claim follows.  
Next, note that the difference in the expectation of $\xi_{\delta}$ and $Lh$ is bounded by a constant that does not depend on $l$ times the occupation time of $(-\delta,\delta)\times(-1/2,1/2)$, which we've already seen goes to zero with $\delta$, uniformly in $l$.  Thus we have that
\[
\E^{(l)}\lb F \int_{s}^{t} \xi_{\delta}(\omega(u))\, du \rb \ra
 \E^{(l)}\lb F \int_{s}^{t} Lh(\omega(u))\, du \rb
\]
as $\delta\searrow 0$, uniformly in $l$.  Because this convergence is uniform in $l$, combining this with Equation~\eqref{Eqn:MoreWeakCon} shows that the second term of Equation~\eqref{Eqn:Martingale} goes to zero as desired.

Now suppose the first possibility in Lemma~\ref{Lemma:SigmaE} holds.  The argument is similar to the previous case.  Without loss of generality, we can assume that $\Gamma_{+}(v_0,v_0)$ is positive at $(x,y)$, since the the roles of $\Gamma_{+}$ and $\Gamma_{-}$ are symmetric.  It follows that $\Gamma(v_0,v_0)$ is uniformly positive ``on one side'' of $\Sigma_0$ in $S'$, where $S'$ is some open neighborhood of $(x,y)$, disjoint from $\{r\leq 1/k_0\}$.  Then $S'\setminus\Sigma_0$ is naturally divided into two disjoint, connected open sets, say $S'_1$ and $S'_2$, and $\Gamma(v_0,v_0)$ is bounded below by a positive constant on one of them, which we can assume is $S'_1$.  Then we can find coordinates $(z_1,\ldots,z_4)$ centered at $(x,y)$ such that $L^{(j(l))}$ and $P^{(j(l),k(l))}$ satisfy the assumptions of Lemma~\ref{Lemma:Occupation} for some choice of constants $R$, $c$, and $C$ independent of $l$, with $d=4$ and $T>t$; and where these coordinates are such that 
\[ 
(0,1)\times (-1,1)^3\subset S'\subset B_R 
\]
and $z_1$ restricted to $S'$ is a constant (non-zero) multiple of the signed distance to $\Sigma_0$, with positive $z_1$ corresponding to $S'_1$ (this implies that $\Sigma_0\cap S'$ is $\{z_1=0\}\cap S'$).  In particular, we note that, in the notation of Lemma~\ref{Lemma:Occupation}, this implies that $a_{1,1}$ is bounded below by a positive constant on $S'_1=\{z_1>0\}\cap S'$.

We can assume that the support of $h$ is contained in $S=(-1/2,1/2)^4$.  In contrast to the previous cases, here the definition of $L$ on $\Sigma_e$ matters.  In particular, we assume that $L|_{\Sigma_0}$ is chosen to be the limit when approached from within $S'_2$ (this is equivalent to having it agree with the optimal orientation-preserving coupling on $\Sigma_0$ if $L$ corresponds to the orientation-preserving coupling on $S'_2$, and similarly for the orientation-reversing case).  We will say more about this assumption below.  This implies that, for each $\delta\in(0,1/2)$, we can find a mollified version of $Lh$, which we again denote $\xi_{\delta}$, with the properties that $\xi_{\delta}$ is continuous, the $\xi_{\delta}$ are bounded uniformly in $\delta$, and $Lh$ and $\xi_{\delta}$ are equal outside of $(0,\delta)\times(-1/2,1/2)^3$.  Weak convergence means that
\[
\E^{(l)}\lb F \int_{s}^{t} \xi_{\delta}(\omega(u))\, du \rb \ra
\E\lb F \int_{s}^{t} \xi_{\delta}(\omega(u))\, du \rb 
\]
as $l\ra\infty$.  In addition, Lemma~\ref{Lemma:Occupation} and our choice of $\xi_{\delta}$ imply that
\[
\E^{(l)}\lb F \int_{s}^{t} \xi_{\delta}(\omega(u))\, du \rb \ra
 \E^{(l)}\lb F \int_{s}^{t} Lh(\omega(u))\, du \rb
\]
as $\delta\searrow 0$, uniformly in $l$.  Because this convergence is uniform in $l$, combining these two equations shows that the second term of Equation~\eqref{Eqn:Martingale} goes to zero as desired.

The only aspect of the argument in the preceding paragraph that needs comment is the possibility of globally defining $L$ on $\Sigma_e$.  However, comparing the arguments in the cases when the first or second possibility of Lemma~\ref{Lemma:SigmaE} holds, we see that the choice of $L$ on $\Sigma_e$ only matters near points $(x,y)\in\Sigma_e$ where the second possibility does not hold for any $H$ and where either $\Gamma_{+}(v_0,v_0)$ or $\Gamma_{-}(v_0,v_0)$ is zero at $(x,y)$ (which, of course, must be true if the second possibility doesn't hold, since otherwise we could just take $H$ to be $\Sigma_0$).  Near any such point $(x,y)$, we know that there exists a neighborhood $S$ of $(x,y)$ such that $\Sigma_0\cap S$ is a smooth hypersurface.  If we choose smooth orthonormal frames for $M$ and $N$ on $S$, then the set of points where $L$ must correspond to the orientation-preserving coupling and the set of points where $L$ must correspond to the orientation-reversing coupling are closed, disjoint sets.  It follows that we can make a global choice of $L$ (on $(M\times N)\setminus\{r=0\}$) such that $L|_{\Sigma_0}$ is what it must be for the above argument to work at all points $(x,y)\in\Sigma_e$ where only the first possibility of  Lemma~\ref{Lemma:SigmaE} holds.  From now on, we will assume that $L$ has been defined on $\Sigma_0$ in a way that satisfies the above description.  As this is all that is needed, we've proved that $P^{k_0}$ is a solution to the martingale problem corresponding to $L$ starting at $(x_0,y_0)$ until $\zeta_{1/k_0}$.

To continue, note that the limit of a convergent sequence $P^{(j(l),k(l))}$ depends only on the tail of the sequence.  Thus we've actually shown that if $j(l)\ra\infty$ and $\liminf k(l)\geq k_0$, then the limit of a convergent subsequence is a solution to the martingale problem corresponding to $L$ until $\zeta_{k_0}$.  Choose $(j(l),k(l))$ so that $j(l)$ and $k(l)$ both go to infinity with $l$, pass to a convergent subsequence, and call the limit $P$.  Then because $\liminf k(l) =\infty$, $P$ is a solution to the martingale problem corresponding to $L$ until $\zeta_{1/k}$ for every positive integer $k$.  Since $\zeta_{1/k}\nearrow \zeta$ as $k\ra\infty$, it follows that $P$ is a solution to the martingale problem corresponding to $L$ until $\zeta$.  The proof of the existence of a solution to the martingale problem for $L$ starting from any point $(x_0,y_0)\in (M\times N)\setminus\{r=0\}$ and stopped at the first hitting time of $\{r=0\}$ is complete.

To prove the final assertion of the theorem, we recall that $L$ is smooth on the complement of $\Sigma_{0}$.  Thus the martingale problem corresponding to $L$ has a unique solution until the first hitting time of $\Sigma_0$.  By uniqueness, the solution just constructed (and, moreover, any solution to the martingale problem for $L$) must agree with this solution until the first hitting time of $\Sigma_0$.  $\Box$

\section{Applications of the coupling}\label{Section:Last}

\subsection{Strong halfspace-type theorems}
A strong halfspace theorem states that two minimal surfaces, satisfying some condition, either intersect or are parallel planes.  Hoffman and Meeks~\cite{HoffmanMeeks} proved a strong halfspace theorem for (complete) properly immersed minimal surfaces.  Their proof used geometric measure theory to show that two such non-intersecting minimal surfaces are separated by a stable minimal surface, which must then be a plane (since planes are the only stable minimal surfaces in $\Re^3$ by a result of Schoen~\cite{Schoen}).  This reduces the problem to the corresponding weak halfspace theorem.  Rosenberg~\cite{Rosenberg} proved a strong halfspace theorem for complete minimal surfaces of bounded curvature, as did Bessa, Jorge, and Oliveira-Filho~\cite{BJO} (this later paper also gives a ``mixed'' strong halfspace theorem in which one minimal surface is properly immersed and the other is complete with bounded curvature).

The ultimate goal of introducing our coupled Brownian motions is to show that the particles couple, either with positive probability or with probability one.  If the particles couple with positive probability, then the minimal surfaces on which they move obviously intersect.  This gives a potential method for proving strong halfspace theorems or similar results.  The issue is proving that the particles couple.  The distance between the particles, under our coupling, is dominated, after time change, by a two-dimensional Bessel process.  Recall that a two-dimensional Bessel process comes arbitrarily close to zero, while a Bessel process of dimension less than two strikes zero in finite time almost surely.  Thus, heuristically, we see that the only obstacles to our particles coupling is that the distance between them might converge, corresponding to the process accumulating only finite quadratic variation, or that the the process might look too much like a two-dimensional Bessel process when the particles are close, causing them to come arbitrarily close but never to couple.  The rate at which the quadratic variation grows and the ratio of the drift to the dispersion are both determined by the configuration of the tangent planes, that is, by the angles $\theta$, $\phi$, and $\psi$ introduced above.

Unfortunately, it's not clear how to get the necessary control of the evolution of these angles to prove the full strong halfspace theorem for either properly immersed or bounded curvature minimal surfaces.  Instead, we have the following partial result.

\begin{THM}\label{THM:WeakWeakHalfspace}
Let $M$ be minimal surface that is either recurrent or stochastically complete with bounded curvature, and let $N$ be a stochastically complete minimal surface.  Then if $M$ is not flat, $\dist(M,N)=0$.
\end{THM}

\emph{Proof:}  Consider Brownian motion on $M$ and $N$, coupled as described above (for arbitrary starting points), and assume that $M$ is not flat.  Then the particles will become arbitrarily close (whether or not they meet) as long as the distance between them does not converge to some positive limit.  To show that this cannot happen, we proceed by contradiction.  Assume that, with positive probability, $r_t$ converges to a positive value (and thus that the process continues for all time).  The vector $x_t-y_t$ is an $\Re^3$-valued martingale, and it is easy to see that its length can only converge if the vector itself converges.  This, in turn, means that the direction of the vector in $\S^2$ must converge.    However, because $M$ is recurrent or stochastically complete with bounded curvature and is not flat, we know from Theorem~\ref{THM:GaussProcess} that, up to a set of probability zero, any path that continues for all time has a normal vector that spends an infinite amount of time in every open set of $\S^2$.  We conclude that, along such a path, the system spends an infinite amount of time where $\theta\in (3\pi/4,\pi/2)$.  The rate of growth of quadratic variation is bounded from below by a positive constant on this set, and thus the $r_t$-process must accumulate infinite quadratic variation.  This contradicts our assumption that it converges, and the proof is complete.  $\Box$

Note that this proof shows that any Brownian motion on $M$ almost surely becomes arbitrarily close to $N$ and vice versa.  Further, in contrast to Rosenberg's result, only one of the minimal surfaces needs to be stochastically complete and have bounded curvature (or be recurrent); the other need only be stochastically complete.  On the other hand, the weakness of this theorem is obvious.  We have only shown that the distance between the surfaces is zero, not that they actually intersect.  As the proof makes clear, the difficulty is controlling the process for small $r$, to rule out both the possibility that $r$ converges to zero and the possibility that it becomes arbitrarily close to zero without ever hitting it.  One might hope that better understanding of the process for small $r$ would allow the theorem to be strengthened to conclude that the surfaces intersect.

In the properly immersed case, we don't even have the analogue of Theorem~\ref{THM:WeakWeakHalfspace}.  This is a consequence of our inability to control the long-term behavior of the normal vector to any extent greater than that implied by the weak halfspace theorem.

\subsection{Maximum principle at infinity}

In a more positive vein, one nice feature of the use of coupled Brownian motions is that one expects results to extend naturally to the case of minimal surfaces with boundary, as mentioned above.  By a minimal surface with boundary, we mean a surface with boundary together with an immersion which is minimal on the interior and extends continuously to the boundary.  Our approach requires extending Theorem~\ref{THM:Main} to the case when one or both of $M$ and $N$ are allowed to have boundary.  To do this, assume that $M$ has non-empty boundary.  Then it is fairly straightforward to show that the first hitting time of the boundary of $M$ is almost surely continuous, with respect to Brownian motion on $M$, on the set where it is finite.  The same is true of $N$ in case it has non-empty boundary.  The first hitting time of the boundary of $M\times N$, which we denote $\eta$, is the minimum of the first hitting times of the boundaries of $M$ and $N$.  All of the processes $P^{(j,k)}$ that we introduced in the proof of Theorem~\ref{THM:Main} and all of their weak limits have marginals that are time-changed (by a constant factor) Brownian motion, so $s\wedge\eta$ and $t\wedge\eta$ (for fixed times $0<s\leq t$) are bounded and almost surely continuous with respect to all of these measures.  Because of this boundedness and almost sure continuity, the weak convergence argument we used in the proof of Theorem~\ref{THM:Main} to show that the limit measure is a solution to the martingale problem for $L$ is compatible with stopping all of the processes at the boundary.  Thus we see that that we can solve the martingale problem for $L$, stopped at the boundary, and this gives an adequate coupling, stopped at the boundary, in the case when one or both of our minimal surfaces has boundary.  (The reason we had to use an approximation argument earlier when stopping the process at $\zeta$ is because it's not clear that $\zeta$ is almost surely continuous.)

The model theorem from (non-stochastic) geometric analysis is the following version of the maximum principle at infinity, recently proved by Meeks and Rosenberg~\cite{MRPreprint}.

\begin{THM}
Let $M$ and $N$ be disjoint, complete, properly immersed minimal surfaces-with-boundary, at least one of which has non-empty boundary.  Then the distance between them satisfies
\[
\dist(M,N)= \min\{\dist(M,\partial N), \dist(\partial M, N)\} .
\]
\end{THM}

This is a generalization of the strong halfspace theorem for properly immersed minimal surfaces.  It is proved by similar methods, although the addition of the boundary makes things more difficult.  We expect the analogue in the bounded curvature case to be true, although to our knowledge no (previous) work has been done in that direction.  However, we have the following version of the maximum principle at infinity for minimal surfaces-with-boundary of bounded curvature.

\begin{THM}\label{THM:MaxAtInfinity}
Let $M$ and $N$ be stochastically complete minimal surfaces-with-boundary, at least one of which has non-empty boundary, such that $\dist(M,N)>0$.  If $M$ has bounded curvature or is recurrent, and is not flat, then
\[
\dist(M,N)= \min\{\dist(M,\partial N), \dist(\partial M, N)\} .
\]
\end{THM}
\emph{Proof:}  Suppose $M$ and $N$ satisfy the hypotheses of the theorem, and that $\dist(M,N)=a>0$.  Consider any point $(x_0,y_0)$ in the interior of $M\times N$. We run a coupled Brownian motion starting at $(x_0,y_0)$, stopped when it hits the boundary.  With one caveat, it is clear from our proof of Theorem~\ref{THM:WeakWeakHalfspace} that $r_t$ almost surely hits the boundary in finite time, having accumulated only finite quadratic variation, since otherwise the process would hit a level below $a$ with positive probability.  The caveat is that the proof of Theorem~\ref{THM:WeakWeakHalfspace} uses Theorem~\ref{THM:GaussProcess}, which we have not proved for surfaces-with-boundary.  However, we can prove that any Brownian path on $M$ with an infinite lifetime has a normal vector which accumulates infinite occupation time in every open set of $\S^2$, up to a set of probability zero, as follows.  If $M$ is recurrent this is clear for the same reasons as before, so assume that $M$ is transient (meaning the interior of $M$ is transient).  Let $\tilde{M}$ be the universal cover of the interior of $M$, so that $\tilde{M}$ is conformally equivalent to the unit disk.  Then $\tilde{M}$ can be described by Weierstrass data just as before, and $\tilde{M}$ also has bounded curvature.  Further, we see that the argument in the proof of Theorem~\ref{THM:GaussProcess} can now be applied to the set of paths with infinite lifetime, so we conclude that any such path has a normal vector which accumulates infinite occupation time in every open set of $\S^2$, up to a set of probability zero.  This establishes our claim that $r_t$ almost surely hits the boundary of $M\times N$ in finite time.

To complete the proof, first suppose that, with positive probability, $r_t$ accumulates no quadratic variation.  It follows that this set of paths produces points $(x,y)$ on the boundary with $\dist(x,y)=\dist(x_0,y_0)$.  For the other case, suppose that $r_t$  almost surely accumulates positive quadratic variation.  Then, by comparison with a two-dimensional Bessel process, there is some $\epsilon>0$ such that $r_t$ hits the level $\dist(x_0,y_0)-\epsilon$ with probability at least $1/2$.  Since the process must stop before $r_t$ gets below level $a$, comparison with a two-dimensional Bessel process shows that there is a positive probability that the paths along which $r_t$ hits $\dist(x_0,y_0)-\epsilon$ are stopped at the boundary before $r_t$ increases to $\dist(x_0,y_0)-\epsilon/2$.  This produces points $(x,y)$ on the boundary with $\dist(x,y)<\dist(x_0,y_0)$.

It follows from the above that, for any point $(x_0,y_0)$ in the interior of $M\times N$, there is a point $(x,y)$ on the boundary of $M\times N$ such that $\dist(x,y)\leq\dist(x_0,y_0)$.  We conclude that $\dist(M,N)= \min\{\dist(M,\partial N), \dist(\partial M, N)\}$, and the theorem is proved.  $\Box$

\subsection{Liouville theorems}

In the previous section, we coupled Brownian motions on two different surfaces in an attempt to control the distance between these surfaces.  However, one can also consider coupling two Brownian motions started at different points on the same surface.  This gives an approach to proving that there are no non-constant bounded harmonic functions on certain classes of minimal surfaces.  In particular, suppose that $M$ is an embedded (by which we mean injectively immersed) minimal surface such that Brownian motions started from any two points couple almost surely.  By embeddedness, the fact that they couple in the extrinsic distance also means that they couple in the intrinsic distance.  Then the standard representation of harmonic functions as integrals of Brownian motion with respect to bounded stopping times shows that any bounded harmonic function must be constant.

Our efforts in this direction are guided by the following conjecture, which appears to go back to Sullivan (see Conjecture 1.6 of~\cite{MeeksSurvey} and the surrounding discussion).

\begin{Conj}
A complete, properly embedded minimal surface admits no non-constant, positive harmonic functions.
\end{Conj}

Though the full conjecture remains open, various special cases are known.  Of course, any class of surfaces which are recurrent satisfies the theorem.  For example, Theorem 3.5 of~\cite{CKMR} states that any complete, properly embedded minimal surface with two limit ends (see the introduction of the paper just cited for a discussion of ends and limit ends) is recurrent.  As for results that apply to transient surfaces, in~\cite{MPR}, Meeks, P\'erez, and Ros prove the conjecture under the additional assumption that the surface possesses one of various symmetries (such as being doubly or triply periodic, or having a sufficiently large isometry group).  Their method is to look at the universal cover, which must be conformally equivalent to the unit disk, and show that the group of Deck transformations must be so large that it can leave only constant positive harmonic functions invariant.

We provide another partial result, under the additional assumption that $M$ has bounded curvature.  Also note that, as indicated above, our result only prohibits non-constant bounded (rather than positive) harmonic functions.

\begin{THM}
Let $M$ be a complete, properly embedded minimal surface of bounded curvature.  Then $M$ has no non-constant bounded harmonic functions.
\end{THM}

\emph{Proof:}  Choose any two distinct points $x_0$ and $y_0$ in $M$, and run our adequate  coupling of Brownian motions from these points.  Our proof of Theorem~\ref{THM:WeakWeakHalfspace} shows that $r_t$ either hits zero in finite time or else spends an infinite amount of time in every neighborhood of zero.  To prove the theorem, it is enough to prove that $r_t$ almost surely hits zero in finite time, as discussed above.  

As a consequence of their maximum principle at infinity (which we stated above), Meeks and Rosenberg were also able to prove that any properly embedded minimal surface with bounded Gauss curvature has a fixed size tubular neighborhood (see the first paragraph of Section 5 of~\cite{MRPreprint} along with Theorem 5.3).  With bounded curvature, the existence of such a tubular neighborhood implies that there is some $a>0$ such that, whenever $r_t=r(x_t,y_t)\leq a$, the distance between $x_t$ and $y_t$ with respect to the metric on $M$, which we denote $\dist_M(x_t,y_t)$, is less than or equal to $2a$.  Further, if this property holds for some particular $a$, then it also holds for any smaller $a$.

Because the curvature is bounded and the embedding is minimal, the entire second fundamental form of the embedding is uniformly bounded.  This means that any fixed sized (with respect to the metric on $M$) neighborhood of a point is uniformly comparable to the tangent plane at that point.  More concretely, for any $\epsilon$, we can choose $a$ small enough so that, whenever $\dist_M(x,y)<2a$, the resulting configuration satisfies
\[\begin{split}
\theta\in\lb\frac{\pi}{2}-\epsilon,\frac{\pi}{2}\rb,\quad \phi\in\lb\frac{\pi}{2}-\epsilon,\frac{\pi}{2}\rb,\quad \psi\in\lb0,\epsilon\rb .
\end{split}\]
One consequence of this is that, by choosing $a$ small enough, we can guarantee that the set of points with $r\leq 2a$ is disjoint from $\Sigma_0$.  In fact, because we also know that the curvature is bounded, the subset of $M\times N$ where $r\leq 2a$ is a positive distance (in the product metric) from $\Sigma_0$.  Thus we can assume that our modified coupling, given by $L$, agrees with the optimal coupling (given by $\ol{L}$) on the set where $r\leq 2a$.  Under the optimal coupling, $\theta=\phi=\pi/2$, $\psi=0$ corresponds to $r_t$ evolving like the standard mirror coupling in the plane, $dr_t = 2\, dB_t$.  So by taking small enough $\epsilon$ and $a$, the configuration can be made arbitrarily close to that of the standard mirror coupling on the plane, which is just time-changed Brownian motion, whenever $\dist_M(x_t,y_t)<2a$.  The other consequence of the fact that the set $\{r\leq 2a\}$ is disjoint from $\Sigma_0$ that we need is that the $(x_t,y_t)$-process, started at any point in $\{ (x,y): r(x,y)\leq 2a\}$, is unique at least until its first exit time from this set.  In particular, the process has the strong Markov property on any  time interval during which it is contained in this set.

The above shows that the rate of growth of quadratic variation of $r_t$ is bounded from below on $\{ (x,y): r(x,y)\leq 2a\}$, so if $r_t$ doesn't strike zero in finite time, it hits the level $a$ and then leaves the set $\{ (x,y): r(x,y)\leq 2a\}$ infinitely many times.  
The point now is to argue that each time $r_t$ hits $a$ it has some probability, bounded from below, of hitting 0 before it hits $2a$.  Then we wish to use the (almost) independence of these events to conclude that, almost surely, on one of these occasions $r_t$ must hit 0.  As indicated above, this proves the theorem.  The rest of the proof is devoted to making this argument precise.

Note that a one-dimensional Bessel process (which corresponds to $\theta=\phi=\pi/2$, $\psi=0$) has the property that, if started at some $l$, it hits zero before it hits $2l$ with probability $1/2$.  
In light of the above, we can choose $a$ so that whenever $r_t\leq 2a$, $r_t$ is dominated by a time-changed Bessel process of dimension $d$, where $d$ is such that a Bessel process of dimension $d$ started at $a$ strikes zero before $2a$ with some probability $p>0$, and so that the time-change satisfies the estimate $d\tau/dt\geq 3/4$.
Let $t_1$ be the first time that $r_t$ hits the level $a$, and let $\rho^{(1)}_{\tau(t)}$ be the comparison Bessel process started from level $a$ at time $t_1$;
\[
\rho^{(1)}_{\tau(t)} = a + W_{\tau(t)}-W_{\tau(t_1)}+\int_{\tau(t_1)}^{\tau(t)} \frac{d-1}{2\rho^{(1)}_{\tau(s)}} \, ds , \quad
\tau(t) = \int_0^t f \, ds ,
\]
where $W_\tau$ is a Brownian motion.

As usual, we stop $(x_y,y_t)$ and $\rho^{(1)}_{\tau(t)}$ if $r_t$ hits zero.  Let $\sigma_1$ be the first time after $t_1$ that $\rho^{(1)}_{\tau(t)}$ hits $2a$; by convention, $\sigma_1=\infty$ if $r_t$ strikes zero before $\rho^{(1)}_{\tau(t)}$ hits $2a$.  We now iterate this procedure.  For $n\geq 2$, let $t_n$ be the first time after $\sigma_{n-1}$ that $r_t$ hits $a$, and let $\sigma_n$ be the first time after $t_n$ that either $r_t$ hits 0 or $\rho^{(n)}_{\tau(t)}$ hits $2a$.  Here $\rho^{(n)}_{\tau(t)}$ is defined by the same equation as $\rho^{(1)}_{\tau(t)}$, except that it is begun at from level $a$ at time $t_n$.  Now consider the processes $\{\rho^{(n)}_{\tau(t)}, t_n\leq t\leq t_{n+1}\}$.  

Recall that, starting from any point with $r=a$, the $(x_t,y_t)$-process is unique at least until $r_t$ hits $2a$.  Using this and the fact that the $[t_n, t_{n+1}]$ are disjoint, we see that this collection of processes is independent, and that each one enjoys the strong Markov property during its lifetime.  Thus, the probability that $\sigma_n$ is finite (that is, the probability that $r_t$ drops to level $a$ and escapes back up to $2a$ without striking zero $n$ times) is less that $(1-p)^n$.  Since this goes to zero as $n\ra\infty$, we see that, almost surely, $r_t$ hits zero in finite time, completing the proof.  $\Box$

\bibliographystyle{amsplain}

\providecommand{\bysame}{\leavevmode\hbox to3em{\hrulefill}\thinspace}
\providecommand{\MR}{\relax\ifhmode\unskip\space\fi MR }
\providecommand{\MRhref}[2]{%
  \href{http://www.ams.org/mathscinet-getitem?mr=#1}{#2}
}
\providecommand{\href}[2]{#2}

\end{document}